\documentclass[11pt, reqno]{amsart}
 \usepackage{amsmath}
 \usepackage{amssymb}
\usepackage{bm}

\DeclareMathAccent{\mathring}{\mathalpha}{operators}{"17}

 \usepackage{color}

\newcommand{\mysection}[1]{\section{#1}
      \setcounter{equation}{0}}

\newcommand\osc{\operatornamewithlimits{osc}}

\newtheorem{theorem}{Theorem}[section]
\newtheorem{lemma}[theorem]{Lemma}

\newtheorem{corollary}[theorem]{Corollary} 

\theoremstyle{definition}
\newtheorem{assumption}{Assumption}[section]
\newtheorem{definition}{Definition}[section]
\theoremstyle{remark}
\newtheorem{remark}{Remark}[section]

\newcommand{\loc}{\text{\rm loc}}

 \makeatletter 
 \def\dashint{%
 \operatorname%
 {\,\,\text{\bf--}\kern-.98em\DOTSI\intop\ilimits@\!\!}}
\def\ninf{\qopname\relax\@empty{inf\phantom{p}\!\!\!}}
\def\nelambda{\dot\lambda}
\def\melambda{\ddot\lambda}

 \makeatother

\newcommand\bbeta{\text{\raise-.2ex\hbox{$\bm{\beta}$}}}

\newcommand\bR{\mathbb{R}}

\newcommand\bS{\mathbb{S}}

\newcommand\cB{\mathcal{B}}

\newcommand\cF{\mathcal{F}}

\newcommand\cN{\mathcal{N}}

\def\sft{{\sf t}}

\newcommand\uR{\underline{R}}

\newcommand\ulambda{\underline\lambda}

\begin{document}

\title[Diffusion processes with drift in $L_{d+1}$]
{On diffusion processes with drift in $L_{d+1}$}

\author{N.V. Krylov}
 
\email{nkrylov@umn.edu}
\address{127 Vincent Hall, University of Minnesota,
 Minneapolis, MN, 55455}
 


\begin{abstract}
This paper is a natural continuation of \cite{Kr_20_2}
and \cite{Kr_21_1}
where strong Markov processes are constructed
in time inhomogeneous setting with Borel measurable
uniformly bounded and uniformly nondegenerate diffusion
and drift in $L_{d+1}(\mathbb{R}^{d+1})$ and some properties
of their Green's functions and probability of passing through narrow
tubes are investigated. On the basis of this  here we study
some further properties of these processes such as
Harnack inequality, H\"older continuity
of potentials, Fanghua Lin estimates and so on.
\end{abstract}

\maketitle

\mysection{Introduction}

Let $\bR^{d}$ be a  Euclidean space of points
$x=(x^{1},...,x^{d})$, $d\geq 2$. Fix some 
 $p_{0},q_{0} \in[1,\infty)$  
such that
\begin{equation}
                                                    \label{5.10.1}
  \frac{d}{p_{0}}+\frac{1}{q_{0}}= 1.
\end{equation}

It is proved in \cite{Kr_20_2} that
It\^o's stochastic equations of the form
\begin{equation}
                                                 \label{11.29.2}
x _{t}=x  +\int_{0}^{t}\sigma (t +s,x_{s})\,dw_{s}
+\int_{0}^{t}b (t +s,x_{s}) \,ds
\end{equation}
admit  weak solutions,
where $w_{s}$ is a $d$-dimensional Wiener process,
$\sigma$ is a uniformly nondegenerate, bounded,
Borel function with values in the set of symmetric
$d\times d$ matrices, $b$ is a Borel measurable $\bR^{d}$-
valued function given on $\bR^{d+1}=(-\infty,\infty)\times\bR^{d}$
such that
\begin{equation}
                                                    \label{5.10.2}
\int_{\bR}\Big(\int_{\bR^{d}}|b(t,x)|^{p_{0}}\,dx\Big)^{q_{0}/p_{0}}\,dt<\infty 
\end{equation}
if $p_{0}\geq q_{0}$ or
\begin{equation}
                                                    \label{5.10.20}
\int_ {\bR^{d}} \Big(\int_{\bR}|b(t,x)|^{q_{0}}\,dt\Big)^{p_{0}/q_{0}}\,dx<\infty 
\end{equation}
if $ p_{0}\leq q_{0}$. 
Observe that the case $p_{0}=q_{0}=d+1$ is not excluded and in this case
the condition becomes $b\in L_{d+1}(\bR^{d+1})$.

The goal of this article is to study the properties of such
solutions or Markov processes whose trajectories
are solutions of \eqref{11.29.2}. In particular,
in Section \ref{section 10.25.2} for more or less general
processes of diffusion type we derive
the power lower estimate for the probability 
to reach at time $T$
a given ball of radius $\rho$
(in case $t=0$, $x=0$ in \eqref{11.29.2}).
This estimate plays a crucial role in proving
the Harnack inequality in Section \ref{section 10.24.1}.
In Section \ref{section 10.25.2} we also prove
that the probability to reach sets of almost full
measure are strictly bigger than zero. This seemingly
weak  statement is also crucial for proving the
Harnack inequality.

Section \ref{section 12.29.2}
is devoted to proving  estimates from below
for the average time spent in space-time sets 
of small measure in terms of a power
of their measure. We also extract some consequences
of these estimates, which help proving the H\"older
continuity of potentials and harmonic functions
in Section \ref{section 10.24.1} and also
allow us to establish the Fanghua Lin estimates
in Section \ref{section 2.17.1}.

Section \ref{section 10.24.1} is devoted to the 
case when our process
is, actually, not just of diffusion type, but a Markov
(time-inhomogeneous) process, whose existence is shown in
\cite{Kr_20_2}. We prove that their resolvent operators
are bounded in $L_{p,q}$. We prove Harnack's inequality
for the caloric functions associated with these processes,
establish that their resolvent operators map $L_{p,q}$
to the set of H\"older continuous functions, and give
some other estimates for the resolvent in the whole space
and in domains. These results extend some of those
in \cite{Ya_20}.

In Section \ref{section 2.17.1} we give some applications
of our results to the theory of linear
parabolic equations. In particular, we prove
the Harnack inequality and H\"older continuity
of their solutions, which in case $p=q=d+1$
are known as Krylov-Safonov estimates and played
an enormous role in the theory of linear
and {\em fully nonlinear\/} elliptic
and parabolic equations with bounded coefficients.
The solutions we are dealing with are of class 
$W^{1,2}_{p,q}$ with $d_{0}/p+1/q\leq 1$ and $d_{0}<d$.
We now have the opportunity to consider
lower order  coefficients in $L_{p,q}$-spaces
and develop $W^{1,2}_{p,q}$-solvability theory.
As is mentioned already, in this section we also
derive the Fanghua Lin estimate,  which is one of  
the starting point of the regularity
 theory
of fully nonlinear equations as presented in \cite{Kr_18}.

The final Section \ref{section 2.14.1} is an appendix
in which we collected some results from
\cite{Kr_21_1} frequently used in the main text.

To the best of the author's knowledge
our results are new even if $p_{0}=q_{0}=d+1$,
with such low integrability of  $b$ and general $\sigma$
H\"older and Harnack properties were unknown.

It is worth mentioning that there is a vast 
literature about stochastic equations when
\eqref{5.10.1} is replaced with $d/p_{0}+2/q_{0}\leq1$.
This condition is much stronger than ours.
Still we refer the reader to the recent articles
\cite{Na_18},
\cite{BFGM_19},  \cite{XZ_20}  and the 
references therein
for the discussion of many powerful 
and exciting results obtained
 under this stronger condition.
There are also many papers when this condition is
considerably relaxed on the account of
imposing various regularity conditions
on $\sigma$ and $b$ and/or considering
random initial conditions with bounded density,
 see, for instance,
\cite{ZZ_19}, \cite{Zh_20} and the 
references therein. Restricting
the situation to the one when $\sigma$ and $b$
are independent of time allows one to
relax the above conditions significantly
further, see, for instance, \cite{KS_19}
and the 
references therein.

 Introduce
$$
B_{R} =\{x\in\bR^{d}:| x|<R\}, 	\quad B_{R}(x)=x+B_{R},
\quad C_{T,R}=[0,T)\times B_{R},
$$
$$
C_{T,R}(t,x)=(t,x)+C_{T,R},\quad C_{R}(t,x)=C_{R^{2},R}(t,x),
\quad C_{R} =C_{R}(0,0),
$$
$$
 D_{i}=\frac{\partial}{\partial x^{i}},
\quad D_{ij}=D_{i}D_{j}\quad \partial_{t}=\frac{\partial}{\partial t}.
$$
For $p,q\in[1,\infty]$ and domains $Q\subset \bR^{d+1}$
 we introduce the space $L_{p,q}(Q)$ as the space of Borel
functions on $Q$ such that
$$
\|f\|^{q}_{L_{p,q}(Q)}:=
\int_{\bR}\Big(\int_{\bR^{d}}I_{Q}(t,x)|f(t,x)|^{p}\,dx\Big)^{q/p}\,dt<\infty 
$$
if $p\geq q$ or
$$
\|f\|^{p}_{L_{p,q}(Q)}:=\int_ {\bR^{d}} \Big(\int_{\bR}I_{Q}(t,x)
|f(t,x)|^{q}\,dt\Big)^{p/q}\,dx<\infty 
$$
if $ p\leq q$ with natural interpretation
of these definitions if $p=\infty$ or $q=\infty$.
If $Q=\bR^{d+1}$, we drop $Q$ in the above notation.
Observe that
$p $ is associated with   $x$ and
$q$ with   $t$ and the interior
integral is always elevated to the power $\leq 1$.
 In case $p=q=d+1$ we abbreviate $L_{d+1,d+1}=L_{d+1}$.
For the set of functions on $\bR^{d}$ summable to the $p$th
power we use the notation $L_{p}(\bR^{d})$.

If $\Gamma$ is a measurable subset of $\bR^{d+1}$ we denote by
$|\Gamma|$ its Lebesgue measure. The same notation
is used for measurable  subsets of $\bR^{d}$ with $d$-dimensional
Lebesgue measure. We hope that it will be clear
from the context which Lebesgue measure is used.
If $\Gamma$ is a measurable subset of $\bR^{d+1}$ and
$f$ is a function on $\Gamma$ we denote 
$$
\dashint_{\Gamma}f\,dxdt=\frac{1}{|\Gamma|}
\int_{\Gamma}f\,dxdt.
$$
In case $f$ is a function on a measurable subset $\Gamma$
of $\bR^{d}$ we set
$$
\dashint_{\Gamma}f\,dx =\frac{1}{|\Gamma|}
\int_{\Gamma}f\,dx .
$$

\mysection{The case of general diffusion type processes
with drift in $L_{p_{0},q_{0}}$}

                                        \label{section 10.25.2}

Let $(\Omega,\cF,P)$ be a complete probability
space, let $\cF_{t}, t\geq0$, be an increasing family of
complete $\sigma$-fields $\cF_{t}\subset\cF$,  
let $w_{t}$ be an $\bR^{d}$-valued Wiener process
relative to $\cF_{t}$. Fix $\delta\in (0,1)$
and denote by $\bS_{\delta}$
the set of $d\times d$ symmetric matrices
whose eigenvalues are between $\delta$ and
$\delta^{-1}$. Assume that we are given
an $\bS_{\delta}$-valued $\cF_{t}$-adapted process
$\sigma_{t}=\sigma_{t}(\omega)$ and an $\bR^{d}$-valued
$\cF_{t}$-adapted process $b_{t}$, such that
$$
\int_{0}^{T}|b_{t}|\,dt<\infty
$$
for any $T\in(0,\infty)$ and $\omega$. Define
$$
x_{t}=\int_{0}^{t}\sigma_{s}\,dw_{s}+\int_{0}^{t}b_{s}\,ds
$$
and for $R\in (0,\infty)$ define $\tau_{R}(x)$
as the first time $(t,x+x_{t})$ exits from $C_{R}$,
$\tau'_{R}(x)$ as the first time $ x_{t} $ exits from $B_{R}$,
$\tau_{R}=\tau_{R}(0)$, $\tau'_{R}=\tau'_{R}(0)$.

\begin{assumption}
                                    \label{assumption 10.20.1}
We are given a function $h\in L_{p_{0},q_{0},\loc}$
such that
$$
|b_{t}|\leq  h(t
,x_{t}).
$$
Furthermore, there exists a  
bounded nondecreasing  function 
$\bar b_{R}$, $R\in(0,\infty)$, such that
 for any $(t,x)\in\bR^{d+1}$ 
    and  $R\in(0,\infty)$ we have
\begin{equation}
                                               \label{8.19.1}
 \|h\|^{q_{0}}_{L_{p_{0},q_{0}}(C_{R}(t,x))} \leq
\bar b_{R}R .
\end{equation} 
\end{assumption}

\begin{assumption}
                                   \label{assumption 12.18.2}
We take $\bar N=\bar N(d,p_{0},\delta)$  introduced
in Theorem \ref{theorem 8.2.1}   and suppose that
there exists $\uR\in(0,\infty)$
such that
\begin{equation}
                                     \label{1.23.2}
\bar N \bar b_{ \uR}< 1.
\end{equation}
\end{assumption}

This assumption as well as Assumption  
\ref{assumption 10.20.1} 
is supposed to hold throughout the article.

\begin{remark}
Throughout the article we fix a number 
$$
\bar R\in[\uR,\infty).
$$ 
In some places we write that certain constants
depend ``only on... $\bar R$ and $\bar b_{\bar R}$''
and the reader might notice that, actually,
the constants depend ``only on... $\bar R$ and $\bar b_{2\bar R}$''.
In such situations it is useful to note that, as is easy to see,
$\bar b_{2\bar R}$ can always be chosen less than or equal
to $N(d)\bar b_{\bar R}$. Also note that if we take
$\bar R=\uR$, then mentioning $\bar b_{\bar R}$ becomes
unnecessary, because $\bar b_{\uR}\leq\bar N^{-1}$.

\end{remark}

Our first big project is to prove a version of
Theorem 4.17 of 
\cite{Kr_19_1},
which provides an important step
toward establishing Harnack's inequality
for caloric functions. It is worth saying that
in the case of bounded $b$ Theorem \ref{theorem 12.7.2}
is proved by constructing a rather simple barrier,
see the PDE argument in the proof of Lemma 9.2.1
(``lemma on an oblique cylinder'')
of \cite{Kr_18} or the probabilistic argument
in the proof of Lemma 2.3 of \cite{Ya_20}.
In our case for the same purpose,
 we needed a rather tedious argument like in
Theorem 4.17 of 
\cite{Kr_19_1}
just to get a good control of the {\em spatial\/}
process $x_{t}$.

Below $\bar\xi=\bar\xi(d,\delta)\in(0,1)$
is taken from Theorem \ref{theorem 8.2.1}.

 \begin{theorem}
             \label{theorem 12.7.2}
Let   $R\in(0,\bar R]$, $ \kappa,\eta\in(0,1)$,
$x,y\in B_{\kappa R}$, 
 and $\eta^{-1} R^{2}\geq t\geq \eta R^{2}$.
Then there exist $N,\nu>0$, depending only on 
$\kappa,\eta,\bar\xi,\uR$, and $\bar R$,
 such that, for any
$\rho\in(0,1]$,
\begin{equation}
                                                 \label{12.9.1}
NP (x+x_{t}\in B_{\rho R}(y),\tau_{   R}'(x)> t )\geq \rho^{\nu }.
\end{equation}
\end{theorem}

The proof of this theorem, given below after 
appropriate preparations,
 follows that of
Theorem 4.17 of 
\cite{Kr_19_1} and, roughly speaking, consists
of splitting the interval $[0,t]$ into several 
parts, estimating the probability that on the first part
the process will reach a neighborhood of $y$ without 
exiting from $B_{R}$,
and then on the consecutive time intervals shrink
the neighborhood with constant coefficient
in such a way as to arrive at time $t$ in $B_{\rho R}(y)$
without exiting from $B_{R}$.

\begin{lemma} 
                    \label{lemma 12.15.1}
Let  $R\in(0,\bar R]$ and
let
  $ \rho_{0}\in(0,1)$, $\theta\in(0,\infty)$, 
and $\kappa\in[1/2,1)$.
Then there exists
$\mu=\mu( \bar\xi, \bar R,\uR,\kappa, \rho_{0},\theta)>0$ such that
\begin{equation}
                                                 \label{12.7.3}
 P(x+x_{\theta R^{2}}\in B_{\rho_{0}\kappa R }(y),
\tau_{ R}(x)> \theta R^{2}  )\geq \mu,
\end{equation}
whenever   
$x,y\in B_{\kappa R}$.
  Furthermore, one can take $\mu>0$
the same if $\bar\xi, \bar R,\kappa$
are fixed and $\rho_{0}$ and $\theta$ vary in 
compact subsets of their ranges.
 
\end{lemma}

Proof. Observe that \eqref{12.7.3} becomes stronger if
$\rho_{0}$ becomes smaller. Therefore we may assume that
\begin{equation}
                               \label{12.7.4}
  \rho_{0}\leq \min\big(\uR/ ( 16\bar R) , \sqrt{\theta / 
(144T_{1} )},\kappa^{-1}-1 \big) \quad (\leq 1/ 16),
\end{equation}
where $T_{1}=T_{1}(\bar\xi)$ is taken from 
Theorem \ref{theorem 1.24.1}.
Then we split the proof into two cases.

{\em Case 1: $|x-y|\leq 3\rho_{0}^{2}\kappa R$\/}.
In that case, owing to $R\leq\bar R$ and $\rho_{0}\leq
\uR/( 6\bar R)$, we have 
$$
|x-y|\leq (1/2)\rho_{0}\kappa (\uR\wedge R).
$$
By Corollary \ref{corollary 1.25.1},
which is applicable since $\rho_{0}\kappa 
(\uR\wedge R)\leq \uR$,
$$
NP(\sup_{t\leq \theta R^{2}}|x+x_{t}-y|<\rho_{0}\kappa 
(\uR\wedge R))
\geq \exp(-4\nu \theta R^{2}/[\rho_{0}\kappa(\uR\wedge R)]^{2}).
$$
The probability here is less than the probability in
\eqref{12.7.3} since $\uR\wedge R \leq R$ and $\rho_{0}
\kappa R\leq (1-\kappa)R$. Furthermore,
$$
   R^{2}/ (\uR\wedge R) ^{2}\leq \bar R^{2}/\uR^{2}
$$
and this proves \eqref{12.7.3} in the first case.

{\em Case 2: $|x-y|\geq 3\rho_{0}^{2}\kappa R$\/}.
Set
$R_{0}=16\rho_{0}^{2}\kappa R$ and note that $|x|+R_{0}
<\kappa R+(1-\kappa)16\rho_{0}R<R$. Similarly,
$|y|+R_{0}<R$. Therefore, the sausage $S_{R_{0}}(x,y)$,
defined as the open convex hull of $B_{R_{0}}(x)\cup B_{R_{0}}(y)$,
 belongs
to $B_{R}$. By Theorem \ref{theorem 1.24.1} with probability
not less than $\pi_{0}^{n}$ before time
$nT_{1}R_{0}^{2}$ the process $x+x_{t}$
will hit $\bar B_{R_{0}/16}(y)$ without exiting
from $S_{R_{0}}(x,y)$, where
$$
n\leq\frac{4|x-y|}{R_{0}}+\frac{1}{4}.
$$
Since $R_{0}<  R$, $|x-y|< 2 R$, and also thanks to
$144 T_{1}\rho_{0}^{2}\leq \theta$, we  have
$$  
nT_{1}R_{0}^{2}\leq T_{1}R_{0}(4|x-y|+ R_{0}/4)
\leq T_{1}R_{0}9   R=144 T_{1}\rho_{0}^{2}\kappa
R^{2}  \leq \theta R^{2}.
$$
By introducing $\gamma$ as the first time $x+x_{t}$
hits $\bar B_{R_{0}/16}(y)$ we conclude that
\begin{equation}
                                               \label{1.27.1}
P(\tau'_{R}(x)>\gamma, \gamma\leq \theta R^{2})\geq \pi_{0}^{n}.
\end{equation}

Observe also that $R_{0}/16=\rho_{0}R_{1}$,
where $R_{1}=\rho_{0}\kappa R\leq \uR$ and at time $\gamma$
the point $x+x_{\gamma}$ is in $\bar B_{\rho_{0}R_{1}}(y)$.
It follows from Corollary \ref{corollary 1.25.1}
that, given that $\gamma<\infty$,
 with probability $\pi_{1}>0$ depending only on $\bar\xi$,
$\rho_{0}$, and $\theta_{1}=\theta R_{1}^{-2}R^{2}$
($\leq \theta\rho_{0}^{-2}\kappa^{-2}$)
the process $x+x_{t}$ will stay in $B_{ R_{1}}(y)$
on the time interval $[\gamma,\gamma+\theta_{1} R_{1}^{2}]$.
Notice that $\gamma+\theta_{1} R_{1}^{2}\geq 
\theta_{1} R_{1}^{2}=\theta R^{2}$.
Along with \eqref{1.27.1} this implies that
$$
P(x+x_{\theta R^{2}}\in B_{ R_{1} }(y),
\tau_{ R}(x)> \theta R^{2}  )\geq \pi_{0}^{n}\pi_{1}>0.
$$
To prove \eqref{12.7.3},
it only remains to recall that $ R_{1}=\rho_{0}\kappa R$.

This proves the first assertion
of the lemma. The second one is obtained by just
inspecting the above proof.
The lemma is proved.

The following is a particular  case of Theorem \ref{theorem 12.7.2}
for $t=\eta R^{2}$.

\begin{lemma} 
                    \label{lemma 12.15.10}
Let $\kappa,\eta\in(0,1)$. 
Then there are  constants $N,\nu>0$, depending only on
$\kappa,\eta,\bar\xi,\uR$, and $\bar R$, such that,  
for any $R\in(0,\bar R],\rho\in (0,1)  $, and $x \in B_{\kappa R}$,
\begin{equation}
                            \label{12.14.2}
NP\big( \tau_{R}(x) > \eta R^{2},
x+x_{\eta R^{2}}\in B_{\rho R} 
\big) \geq \rho^{\nu }.
\end{equation}
\end{lemma}

Proof. 
We may assume that $\kappa\in(1/2,1)$. 
  Let $\rho_{0}$  
be a   positive solution of
$$
\rho_{0}=\min\Big(\frac{\uR}{16\bar R},
\frac{\sqrt{\eta(1-\rho_{0})}}{12\sqrt{T_{1}}},\frac{1}{\kappa}-1\Big),
$$
observe that it suffices to prove
\eqref{12.14.2} for $\rho\leq\kappa $,
and  set
$$
 n(\rho)=\Big\lfloor
\frac{\ln(\rho/\kappa)}{\ln \rho_{0} }\Big\rfloor+1\quad(\geq1),\quad \bar\eta
=\eta\frac{1-\rho_{0}}{1-\rho_{0}^{2n(\rho)}}.
$$

Note that $\bar\eta\in [\eta(1-\rho_{0}),\eta
(1+\rho_{0})^{-1}]$ and
$$
\rho_{0}\leq \min\Big(\frac{\uR}{16\bar R},
\frac{\sqrt{\bar\eta}}{12\sqrt{T_{1}}},\frac{1}{\kappa}-1\Big)
$$
so that by Lemma \ref{lemma 12.15.1} 
 estimate \eqref{12.7.3},
with $\theta=\bar\eta$, $y=0$, 
is valid and means that  
\begin{equation} 
                         \label{12.16.10}
 P\big(x+x_{ \bar\eta  R^{2}}\in B_{\rho_{0}\kappa R },
\sup_{s\leq  \bar\eta  R^{2}}|x+x_{s}|<   R  \big)\geq \mu,
\end{equation}
whenever $R\in(0,\bar R]$ and  
$x\in B_{\kappa R }$. For $n=1,2,...$ 
introduce ($t_{0}:=0$)
$$
R_{n}=\rho_{0}^{n-1}R ,\quad s_{n}= 
 \bar\eta  R_{n}^{2}=
 \bar\eta R^{2} \rho_{0}^{2(n-1)},\quad
t_{n}=\sum_{k=1}^{n}s_{k},
$$
$$
A_{n}=\{\sup_{s\leq s_{n}}|x+x_{s+t_{n-1}}|< R_{n}\},\quad
\Pi_{n}=\bigcap_{k=1}^{n}A_{k} 
$$
and observe that by the conditional  version of \eqref{12.16.10}
on the set $\{y:=x+x_{t_{n-1}}\in B_{\kappa R_{n}}\}$ we have (a.s.)
\begin{equation}
                                                 \label{12.9.6}
P\Big(y+(x_{t_{n}}-x_{t_{n-1}})\in B_{
\kappa R_{n+1}},
\sup_{s\leq s_{n}}|y+(x_{t_{n-1}+s}-x_{t_{n-1}})|<   R_{n}
\mid \cF_{t_{n-1}}\Big)\geq \mu.
\end{equation}
Furthermore, obviously, for $n\geq 2$,
$$
P^{n}:= 
P(x+x_{t_{n}}\in B_{\kappa R_{n+1 }},\Pi_{n}  )
$$
$$
\geq P(x+x_{t_{n-1}}\in B_{\kappa R_{{\color{blue}n}}},\Pi_{n-1},
$$
$$
x+x_{t_{n-1}}+(x_{t_{n}}-x_{t_{n-1}})\in B_{\kappa R_{n+1}},
\sup_{s\leq s_{n}}|x+x_{t_{n-1}}+(x_{t_{n-1}+s}-x_{t_{n-1}})|<  R_{n} ),
$$
which in light of \eqref{12.9.6} yields $P^{n}\geq \mu P^{n-1}$
and, since for $|x|<\kappa R$
we have $P^{1}\geq \mu$ by \eqref{12.16.10}, it holds that
for $|x|<\kappa R$ and all 
$n\geq   1$ 
\begin{equation}
                                                 \label{12.16.1}
P\big(x+x_{t_{n}}\in B_{\kappa R_{n+1 }},\sup_{s\leq t_{n}}|x+x_{s}|
< R  \big)\geq \mu^{n }.
\end{equation}

Observe that 
$$
t_{n(\rho)}=\bar\eta R^{2}\frac{1-\rho_{0}^{2n(\rho)}}
{1-\rho_{0} }=\eta R^{2},\quad
\kappa R_{n(\rho)+1}= \kappa\rho_{0}R
^{n(\rho)}\leq \rho R.
$$
Therefore, \eqref{12.16.1} implies that
$$
P(\tau_{1}(x)>\eta R^{2},|x+x_{\eta R^{2}}|\leq \rho R)
\geq \mu^{n(\rho) }.
$$
Now to finish the proof,
it only remains to note that
$$
\mu^{n(\rho) }\geq
\mu \exp\Big(
\frac{\ln(\rho/\kappa)}{\ln \rho_{0} }\ln\mu\Big)
=N\rho^{\nu}.
$$
The lemma is proved.

{\bf Proof of Theorem \ref{theorem 12.7.2}}.
Let $R_{1}=(1-\kappa)R$ and note that 
 $\xi:= t/R_{1}^{2}-\eta$ satisfies
$$
\eta^{-1}(1-\kappa)^{-2}>\xi\geq \eta\big[(1-\kappa)^{-2}-1\big].
$$

By the conditional  version of
Lemma \ref{lemma 12.15.10} on the set
$\{z:= x+x_{\xi R_{1}^{2}}\in B_{\kappa R_{1}}(y)\}$ we have (a.s.)
$$
NP\Big(\sup_{s\in [\xi R_{1}^{2},\xi R_{1}^{2}+\eta R_{1}^{2}]}
|z+x_{s}-x_{\xi R_{1}^{2}}-y|<R_{1},
$$
$$
x+x_{\xi R_{1}^{2}+\eta R^{2}_{1}}\in B_{\rho R_{1}}(y)
\mid \cF_{\xi R_{1}^{2}}\Big)
\geq \rho^{\nu}.
$$

By Lemma \ref{lemma 12.15.1}, where we take $\rho_{0}=R_{1}/R$
and replace $\theta$ there with $\xi(1-\kappa)^{2}$,
$$
P(\sup_{s\leq \xi R_{1}^{2}}|x+x_{s}|< R,
x+x_{\xi R_{1}^{2}}\in B_{ \kappa R_{1} }(y) 
 )\geq \mu.
$$
 
By combining these two facts
and using   that
$\xi R^{2}_{1}+\eta R_{1}^{2}=t$, we obviously come
to \eqref{12.9.1}. The theorem is proved.

 Theorem \ref{theorem 11.8.1}
originated in \cite{KS_79} in case $b$ is bounded
and is one of two most important ingredients in the proof
of the Harnack inequality. Observe that in this theorem
we do not claim that
$q(\xi)\ne 0$ for $\xi$ not close to one. This fact
will be proved later.   

\begin{theorem}
                                     \label{theorem 11.8.1}
For any $\kappa\in(0,1)$
there is a 
function $q(\xi)$, $\xi\in(0,1)$,
depending only on $\kappa,\delta,d,\uR, p_{0},\bar R$, 
and, naturally, on $\xi$,
such that for any $R\leq \bar R$, $x\in B_{\kappa R}$,
and closed $\Gamma\subset \bar C_{ R}$ satisfying
$|\Gamma|\geq \xi|C_{ R}|$ we have
\begin{equation}
                               \label{1.7.2}
P (\tau_{\Gamma}(x) <  \tau_{   R}(x) )\geq q(\xi),
\end{equation}
where $\tau_{\Gamma}(x) $ is the first time the process 
$(t, x+x_{t})$
hits $\Gamma$  (and  
 $\tau_{R}(x) $ is its first exit time from
$C_{ R}$). Furthermore, 
$q(\xi)\to 1$ as $\xi\uparrow 1$.
Finally, for any closed $\Gamma'\subset B_{ R}$ satisfying
$|\Gamma'|\geq \xi|B_{ R}|$ we have
\begin{equation}
                               \label{1.7.1}
P (\tau'_{\Gamma}(x) <  \tau'_{   R}(x) )\geq q(\xi),
\end{equation}
where $\tau'_{\Gamma}(x) $ is the first time the process 
$ x+x_{t} $
hits $\Gamma$  (and  
 $\tau'_{ R}(x) $ is its first exit time from
$B_{ R}$).
\end{theorem}

Proof.
By considering what is going on in $B_{(1-\kappa)R}(x)$
we convince ourselves that we may assume that $x=0$.
Also, obviously we may assume that $R\leq\uR$.

 In that case for any $\varepsilon\in(0,\bar \xi/4)$ we have
($\tau_{\Gamma}=\tau_{\Gamma}(0)$)
$$
P(\tau_{\Gamma} > \tau_{R} )\leq
P\Big(\tau_{R} =\int_{0}^{\tau_{R} }I_{C_{R}
\setminus\Gamma}(t,x_{t})
\,dt\Big)
$$
$$
\leq P(\tau_{R} \leq\varepsilon R^{2} )+\varepsilon^{-1}R^{-2} 
E\int_{0}^{\tau_{R}}I_{C_{R}\setminus\Gamma}(t,x_{t}) 
\,dt.
$$
In light of Theorems \ref{theorem 8.20.1} and 
Lemma \ref{lemma 8.16.1}
we can estimate the right-hand side and then obtain
$$
P(\tau_{\Gamma}> \tau_{R})\leq Ne^{-1/(N\varepsilon)}
+N\varepsilon^{-1}R^{d/(d+1)-2}|C_{R}\setminus\Gamma|^{1/(d+1)} 
$$
$$
\leq Ne^{-1/(N\varepsilon)}   
+N\varepsilon^{-1}(1-\xi)^{1/(d+1)}
$$
where the constants $N$  depend only on $d,\delta,p_{0}$.
By denoting
$$
q(\xi)=1-\inf_{\varepsilon\in(0,\bar\xi/4)}\big(
Ne^{-1/(N\varepsilon)}
+N\varepsilon^{-1}(1-\xi)^{1/(d+1)}\big),
$$
we get what we claimed about \eqref{1.7.2}.

Estimate \eqref{1.7.1} follows from \eqref{1.7.2}
if one takes in the latter $\Gamma=[0,R^{2}]\times
\Gamma'$ and observes that
$$
\{\tau _{\Gamma}(x) <  \tau _{   R}(x)\}
\subset \{\tau'_{\Gamma}(x) <  \tau'_{   R}(x)\}.
$$
The lemma is proved.

\mysection{Estimating time spent in space-time sets 
of small measure}
                                  \label{section 12.29.2}

The central result of this section is Theorem \ref{theorem 12.21.1}
which needs some auxiliary constructions and assertions.

We present 
extensions to the case that $b\in L_{d+1}$  of
 probabilistic versions of some
PDE results found in \cite{KS_80}, \cite{Sa_80},  
\cite{Kr_18}. Recall the notation introduced in
the Introduction and also
introduce
$$
C^{o}_{T,R}=(0,T)\times B_{R},\quad C^{o}_{T,R}(t,x)=(t,x)+C^{o}_{T,R},
\quad C^{o}_{R}(t,x)=C^{o}_{R^{2},R}(t,x),
  $$
 $C^{o}_{R}=C^{o}_{R}(0,0)$.
Fix 
 $$
q, \eta ,  \kappa\in(0,1).
$$

 For cylinders $Q= C^{o}_{ \rho}(t ,x ) $
 define       
$$
Q'=(t ,x )-
C^{o}_{\eta^{-1} \rho^{2},\rho},\quad Q''=\big
(t -\eta^{-1}\rho^{2},x \big)
+C^{o}_{\eta^{-1}\rho^{2}\kappa^{2},\rho\kappa},
 $$
$$
Q'_{+}=Q\cup Q'\cup \big(\{t\}\times B_{\rho}(x)\big).
$$

Imagine that the $t$-axis is pointed up vertically. Then 
 $Q'$ is   adjacent to $Q$ from below, the 
two cylinders have a common base, and along the $t$-axis
$Q'$ is $\eta^{-1}$ times longer than $Q$.  
 The cylinder
$Q''$ is obtained by  contracting $Q'$  
to    the
center  of  its  lower base   with the contraction factor  
$\kappa^{-2}$ for the $t$-axis and $\kappa^{-1}$ for the spatial
axes. 

\begin{remark}
                                                   \label{remark 2,14,1}
If  $Q= C^{o}_{ \rho}(t ,x )$, then
  the   distance
between $Q$ and $Q''$ along the $t$ axis is     
\begin{equation}
                                                          \label{2,14,1}
\eta^{-1}\rho^{2}-\eta^{-1}\rho^{2}\kappa^{2}=\eta^{-1}\rho^{2}
(1-\kappa^{2}),
\end{equation}
 which is 1 if $\eta=1-\kappa^{2}$.

\end{remark}

 Let $\Gamma$ be a measurable subset of $C_{1}$
and
 introduce   $\mathcal{B}=\cB(\Gamma,q)$ 
 as the family  of   
{\em open\/}
cylinders $Q$ of  type 
$ C^{o}_{ \rho}(t_{0},x_{0})$ such that   
$$
Q\subset  C_{1 } \quad\text{and}\quad
|Q\cap\Gamma|\ge q|Q|.
$$ 

Finally,  define 
$$
\Gamma''=\bigcup_{Q\in\mathcal{B}}Q'',\quad
\Gamma''_{\varepsilon}=\bigcup_{Q\in\mathcal{B}:|Q|\geq
\varepsilon}Q''.
 $$

 Observe that for $Q\in\mathcal{B}$ the set 
   $Q''$ is open. Hence, 
$\Gamma''$ is open and measurable. 
\begin{lemma}
                                                \label{lem:4.1.6} 
If $|\Gamma|\le q|C_{1 }|$, then   
$$|\Gamma''|\ge\Big(1-\frac{1-q}{3^{d+1}}\Big)^{-1}
(1+\eta)^{-1}\kappa^{d+2}|\Gamma|
 $$
and, if the factor of $|\Gamma|$ above
is strictly bigger than one,
 there exists $\theta=\theta(d,q,\eta,\kappa)>1$
such that for any sufficiently small $\varepsilon>0$
there exists a closed $\Gamma_{\varepsilon}
\subset \Gamma''_{\varepsilon}$ such that
\begin{equation}
                                                 \label{12.21.3}
|\Gamma_{\varepsilon}|\geq\theta |\Gamma|.
\end{equation}
\end{lemma}

The first assertion of the lemma originated
in \cite{KS_80}, \cite{Sa_80},  is presented, for instance
as Lemma 9.3.6 in
\cite{Kr_18}. The second one is proved
in the same way as the second assertion
of Lemma 4.8 of \cite{Kr_19_1}.

\begin{lemma}
                                         \label{lemma 12.20.2}
There is a constant 
$q_{0}=q_{0}(d,\delta,p_{0},\uR,\bar R)
\in(0,1)$ such that for any $R\leq \bar R$,
 Borel set $\Gamma\subset C_{R} $
satisfying $|\Gamma|\geq q_{0}|C_{R}|$,
and $x\in \bar B_{R/2}$ we have
\begin{equation}
                                                    \label{12.20.3}
E\int_{0}^{\tau_{R}(x)}I_{\Gamma}(t ,x+x_{t})
\,dt\geq \mu_{0} R^{2} ,
\end{equation}
where   $\mu_{0}=\mu_{0}(d,\delta,p_{0},\bar R,\uR,
\bar b_{\bar R})\in(0,1)$.

\end{lemma}

Proof.  Note that
in light of
 Corollary \ref{corollary 7.29.1} we have
$   E \tau_{R}(x)\geq \nu  R^{2}$, where $\nu
=\nu(d,\delta,\bar R,\uR )>0$.
By using Theorem \ref{theorem 9.5.1} we get that
$$
E \tau_{R}(x)-
E\int_{0}^{\tau_{R}(x)}I_{\Gamma}(t,x+x_{t})
\,dt=E \int_{0}^{\tau_{R}(x)}I_{C_{R}\setminus
\Gamma}(t,x+x_{t})
\,dt
$$
$$
\leq NR^{(2d_{0}-d)/(d_{0}+1)}(|C_{R}|-|\Gamma|)^{1/(d_{0}+1)}
$$
$$
\leq NR^{2}(1-q_{0})^{1/(d_{0}+1)}\leq N(1-q_{0})^{1/(d_{0}+1)}
E \tau_{R}(x),
$$
where the constants $N$ depend only on $ d,\delta$, $\bar R$,
 $\uR$, $p_{0}$,
and $\bar b_{\bar R}$ and $d_{0}=
d_{0}(d,\delta,\uR,p_{0} )\in(d/2,d)$ is taken from 
\cite{Kr_21_1}. We see how to choose
$q_{0}$ to get the desired result.
The lemma is proved.  

In Lemma \ref{lemma 12.20.3} by $q_{0}$ we mean the one from Lemma \ref{lemma 12.20.2}.

\begin{lemma}
                                         \label{lemma 12.20.3}
Take   $Q= C^{o}_{ \rho}(s,y)$
with $\rho\leq\bar R$,
use the notation $Q',Q'',Q'_{+}$ introduced above,
assume than $\eta=1-\kappa^{2}$,
and suppose that Borel $\Gamma\subset Q$ is such that
$|\Gamma|\geq q_{0}|Q|$. Then there is a constant
$\nu_{0}>0$, depending only on $ \kappa$,
$ d,\delta$, $\bar R$,
 $\uR$, $p_{0}$,
and $\bar b_{\bar R}$,
such that for any $(t_{0},x_{0})\in Q''$
\begin{equation}
                                                    \label{12.20.4}
E \int_{0}^{\tau }I_{\Gamma}(t_{0}+t,x_{0}+x_{t})
\,dt\geq \nu_{0}E \tau,
\end{equation}
where $\tau$ is the first exit time of $(t_{0}+t,x_{0}+x_{t})$
from $Q'_{+}$.
\end{lemma}

Proof. Thanks to   Remark \ref{remark 2,14,1}
we have $s-t_{0}\in(\rho^{2}, \eta^{-1}\rho^{2})$.
Also $|y-x_{0}|<\kappa\rho$. It follows by
 Theorem \ref{theorem 12.7.2} that
$$
P \big(\sup_{r\in[0,s-t_{0}]}|x_{0}+x_{r}-y|< \rho,
|x_{0}+x_{s-t_{0}}-y|< (1/2)\rho\big)\geq \nu ,
$$
where $\nu=\nu(\kappa,d,\delta,\uR,\bar R)>0$.

Next,
for $\gamma$ defined as the first exit time of $(t_{0}+t,x_{0}+x_{t})$
from $Q'$ in light of Lemma \ref{lemma 12.20.2} we have
$$
E \int_{0}^{\tau }I_{\Gamma}(t_{0}+t,x_{0}+x_{t})
\,dt=E I_{\gamma>s-t_{0}}
\int_{\gamma}^{\tau }I_{\Gamma}(t_{0}+t,x_{0}+x_{t})
\,dt
$$ 
$$
\geq EI_{\gamma>s-t_{0},|x_{0}+x_{s-t_{0}}-y|< \rho/2}
E\Big(\int_{0}^{\tau }I_{\Gamma}(t_{0}+t,x_{0}+x_{t})
\,dt\mid \cF_{s-t_{0}}\Big)
$$
$$
\geq \mu_{0}\rho^{2}
P \big(\sup_{r\in[0,s-t_{0}]}|x_{0}+x_{r}-y|< \rho,
|x_{0}+x_{s-t_{0}}-y|< \rho/2\big)\geq \mu_{0}\nu\rho^{2}.
$$
On the other hand, the height of $Q'_{+}$ is $(1+\eta^{-1})
\rho^{2}$, so that $(t_{0}+t,x_{0}+x_{t})$ cannot spend
in $Q'_{+}$ more time than $(1+\eta^{-1})
\rho^{2}$. This proves the lemma.  

\begin{lemma}
                                             \label{lemma 2.11.1}
Denote
$$
G_{R}(\Gamma,x):= E\int_{0}^{\tau'_{R}(x)\wedge(2R^{2})}
I_{\Gamma}(t,x+x_{t})\,dt,
$$
fix $q,\kappa\in(0,1)$,
and introduce $\mu_{R}(q)$  as $R^{-2}$ times
the infimum of $G_{R}(\Gamma,x)$
over all Borel $\Gamma\subset C_{R}(R^{2},0)$ satisfying $|\Gamma|
\geq q|C_{R}(R^{2},0)|$  over all $x\in B_{\kappa R}$,
and over all processes $x_{t}$ satisfying our assumptions
with the same $\delta$ and $\bar b_{R}$. Then $\mu_{R}(q)$
is a decreasing function of $R$.
\end{lemma}

The proof  of this lemma, left  to the reader, is easily achieved
by using the self-similar dilations: $x_{t}\to cx_{t/c^{2}}$,
which preserves (actually, makes smaller) $\bar b_{R}$
(see, for instance, the proof of our Theorem \ref{theorem 2.3.1}
in \cite{Kr_21_1}).

\begin{theorem}
                                             \label{theorem 12.21.1}
For
 any $\kappa\in(0,1)$ there exist  $\gamma\in(0,1)$
and $N$,
depending only on $\kappa,d,\delta,p_{0},\uR$,
with $N$ also depending on $\bar R$ and
and $\bar b_{\bar R}$, such that
for any $R\in(0,\bar R]$, $q\in(0,1)$,   Borel
$\Gamma\subset C_{R}(R^{2},0)$ satisfying
$|\Gamma|\geq q |C_{R}(R^{2},0)|$, and $x\in B_{\kappa R}$
we have
\begin{equation}   
                                                   \label{10.1.10}
  E\int_{0}^{\tau'_{R}(x)\wedge(2R^{2})}
I_{\Gamma}(t,x+x_{t})\,dt\geq N^{-1}q^{1/\gamma}R^{2}.
 \end{equation}

\end{theorem}

Proof. By using the notation from Lemma \ref{lemma 2.11.1},
our assertion is rewritten as
\begin{equation}   
                                                   \label{2.12.1}
\mu_{R}(q)\geq N^{-1}q^{1/\gamma}R^{2}.
\end{equation}
 
Fix $\kappa\in(0,1)$,  perhaps larger than the one in the statement
of the theorem, such that for $\eta=1-\kappa^{2}$
and $q=q_{0}$ the factor of $|\Gamma|$ in 
Lemma \ref{lem:4.1.6} is strictly bigger than one and
take $\theta=\theta(d,q_{0},1-\kappa^{2},\kappa)>1$
from that lemma. 

Next, observe that a combination of Lemma \ref{lemma 12.20.2}
and Theorem \ref{theorem 12.7.2}, as in the proof of Lemma
\ref{lemma 12.20.3}, leads to the conclusion that
there exists $ \mu_{0}\in(0,1)$, depending only on
$ \kappa,d,\delta,p_{0},\uR,\bar R$, and $\bar b_{\bar R}$,
 such that
$$    
\mu_{R}(q)\geq \mu_{0} 
$$
for $q\in [q_{0},1]$ and $R\leq\bar R$.

We will be comparing $\mu_{R}(q')$ and $\mu_{R}(q'')$
for $ 0 <q'<q''<1$ such that
\begin{equation}
                                               \label{1.5.1}
(1+\theta) q'\geq 2q''.
\end{equation}

We take a Borel
$\Gamma\subset C_{R}(R^{2},0)$ satisfying
$|\Gamma|\geq q' |C_{R}(R^{2},0)|$ and
in the construction before Lemma \ref{lem:4.1.6} we replace $C_{1}$
by $C_{R}(R^{2},0)$, keep all other notation, and
from the chosen $\Gamma,\kappa,\eta$, and $q_{0}$ (not $q'$) we  
 build up  the
set  $\Gamma _{\varepsilon}$ and take
$\varepsilon$ so small that \eqref{12.21.3} holds.
  There are  two cases:  

(i) $\big|\Gamma_{\varepsilon}\setminus 
C_{R}(R^{2},0)\big|\leq (q''-q')|C_{R}(R^{2},0)|$,

(ii) $\big|\Gamma_{\varepsilon}\setminus
C_{R}(R^{2},0)\big|> (q''-q')|C_{R}(R^{2},0)|$.

{\em Case  (i )\/}. Our  goal is to show that
\begin{equation}
                                                 \label{7,26,5}
G_{R} (\Gamma,x)
\geq\min\big(\mu_{0}, \nu_{0}\mu_{R}(q'')\big)R^{2},\quad|x|\leq \kappa R,
\end{equation}
where $\nu_{0}$
  depends only on  $\kappa, d,p_{0},\delta,\uR,\bar R$,
and $\bar b_{\bar R}$.

Observe that, if $|\Gamma|\geq q_{0}|C_{R }|$, by definition
$G _{R}(\Gamma,x)\geq \mu(q_{0})R^{2}\geq\mu_{0}R^{2}$ 
for $|x|\leq \kappa R$.
 Hence, we may assume that
$$
|\Gamma|< q_{0}|C_{1 }|.
$$
In that case define
$$
\hat\Gamma_{\varepsilon}=\Gamma_{\varepsilon}\cap C_{R }(R^{2},0).
 $$

Notice that by definition and  Lemma \ref{lem:4.1.6}   
$$
q'|C_{R }|\leq|\Gamma|\leq \theta^{-1}|\Gamma_{\varepsilon}|.
 $$
Moreover, by assumption   
$$
|\Gamma_{\varepsilon}|=\big|\Gamma_{\varepsilon}
\setminus C_{R }(R^{2},0)\big|+|\hat  \Gamma_{\varepsilon}|
\leq (q''-q')|C_{R }|+|\hat  \Gamma_{\varepsilon}|.
 $$
Due to \eqref{1.5.1}, it follows that
$$
|\hat \Gamma_{\varepsilon}|\geq q''|C_{R }|,
$$
so that      
$$
G_{R}(\hat \Gamma_{\varepsilon},x)\geq\mu_{R}(q'')R^{2} ,
\quad|x|\leq \kappa R.
 $$

We now estimate $G_{R}( \Gamma ,x)$ from
below by means of $G_{R}(\hat \Gamma_{\varepsilon},x)$
using Lemma \ref{lemma 12.20.3}. Since $\Gamma_{\varepsilon}
\subset \Gamma''_{\varepsilon}$,
the closed set $\Gamma_{\varepsilon}$ is 
covered by the family $\{Q'':Q\in\mathcal{B},|Q|\geq
\varepsilon\}$. Then there is
finitely many $Q(1),...,Q(n)\in \mathcal{B}$ such that
$|Q(i)|\geq\varepsilon$, $i=1,...,n$, and
$$
\Gamma_{\varepsilon}
\subset \bigcup_{i=1}^{n}Q''(i)
=:\Pi_{\varepsilon}.
$$

Then for $(t,x)\in \Pi_{\varepsilon}$ define  
$i (t,x) $
as the first $i\in\{1,...,n\}$ for which
$(t,x)\in Q''(i)$. 
Also set $Q'_{+}(0)=C_{2R^{2},R}$
and $i(t,x)=0$ if $(t,x)\in\partial C_{2R^{2},R}$.
Now define recursively $\gamma^{0}=0$,
$\tau^{1}$ as the first time after $\gamma^{0}$ when $(t,x+x_{t})$ exits
from $C _{2R^{2},R}\setminus \Gamma _{\varepsilon}$,
$\gamma^{1}$ as the first  time after $\tau^{1}$
when $(t,x+x_{t})$ exits from $Q' _{+}(i(\tau^{1},x+x_{\tau^{1}}))$,
and generally, for $k=2,3,...$ define
$\tau^{k}$ as the first time after $\gamma^{k-1}$ when $(t,x+x_{t})$ exits
from $C_{2R^{2},R}\setminus \Gamma _{\varepsilon}$,
$\gamma^{k}$ as the first  time after $\tau^{k}$
when $(t,x_{t})$ exits from $Q'_{+}(i(\tau^{k},x+x_{\tau^{k}}))$.
It is easy to check that so defined
$\tau^{k}$ and $\gamma^{k}$ are stopping times
and, since $|Q(i)|\geq\varepsilon$ and the trajectories 
of $(t,x+x_{t})$ are continuous,
$\tau^{k}\uparrow \tau'_{R}(x)\wedge 2R^{2}$ as $k\to\infty$.
Furthermore,
(a.s.) all the $\tau^{k}$'s equal $\tau'_{R}(x)\wedge 2R^{2}$
 for all large $k$.

For a domain $Q\subset\bR^{d+1}$
we denote by $\gamma(s,y,Q)$ the first exit time
of $(s+t,y+x_{s+t}-x_{s})$ from $Q$ and   obtain
$$
G_{R}(\Gamma,x)\geq \sum_{k=1}^{\infty}E\int_{\tau^{k}}
^{\gamma^{k}}I_{\Gamma}(t,x+x_{t})\,dt
$$
$$
=\sum_{k=1}^{\infty}EE\Big(\int_{0}
^{\gamma(s,y,Q_{+}'(i))}I_{\Gamma}(s+t,x_{s+t}-x_{s}+y)\,dt
\mid \cF_{s}\Big)\Big|_{
i=i(s,y),y=x+x_{s},s=\tau^{k}}.
$$
We estimate the interior expectation from below
by Lemma \ref{lemma 12.20.3} and get that $G_{R}(\Gamma,x)/\nu_{0}$
is greater than or equal to
$$
\sum_{k=1}^{\infty}E E\Big(\int_{0}
^{\gamma(s,y,Q'_{+}(i))}I_{\Pi_{\varepsilon}}
(s+t,x_{s+t}-x_{s}+y)\,dt
\mid \cF_{s}\Big)\Big|_{
i=i(s,y),y=x+x_{s},s=\tau^{k}}
$$
$$
\geq  
\sum_{k=1}^{\infty}E E\Big(\int_{0}
^{\gamma(s,y,Q'_{+}(i))} I_{\Gamma_{\varepsilon}}(s+t,x_{s+t}-x_{s}+y)\,dt
\mid \cF_{s}\Big)\Big|_{
i=i(s,y),y=x+x_{s},s=\tau^{k}}
$$
$$
=  
\sum_{k=1}^{\infty}E  \int_{\tau^{k}}
^{\gamma^{k}}I_{\Gamma_{\varepsilon}}( t,x+x_{t})\,dt=
 G_{R}(\Gamma_{\varepsilon},x)\geq
 G_{R}(\hat \Gamma_{\varepsilon},x)\geq \mu_{R}(q'')R^{2}.
$$
This proves \eqref{7,26,5}.

{\em Case (ii)\/}. Here the goal is to prove that  
\begin{equation}
                                                   \label{7,27,1}
G_{R}(\Gamma,x)\geq \mu_{0}\nu\eta^{n}(q''-q')^{n}R^{2},
\quad |x|\leq \kappa,
 \end{equation}
  where $\nu>0$ and $n\geq1$ depend 
only on $d,\delta,\bar R$, $\uR$, and $\kappa$.

First we claim that for some $(t,x)\in\Gamma_{\varepsilon}$
it holds that $t<(q'-q''+1)R^{2}$. Indeed, otherwise
$$
\Gamma_{\varepsilon}\setminus C_{R}(R^{2},0)
\subset C_{(q''-q')R^{2},R}((q'-q''+1)R^{2},0)
$$
 and 
$|\Gamma_{\varepsilon}\setminus C_{R}(R^{2},0)|\leq
(q''-q')|C_{R}|$. It follows that there is a cylinder
$$
Q=C^{o}_{ \rho}(s,y)\in\cB
$$
 such that $Q'$
contains points in the half-space $t<(q'-q''+1)R^{2}$. 
Since $q'<q''$, we have $q'-q''+1<1$, and since $Q'$ is adjacent
to $Q\subset C_{R }(R^{2},0)$, this implies that the height of $Q'$
is at least $(q''-q')R^{2}$, that is,  
\begin{equation}
                                            \label{12.22.7}
\rho^{2}\eta^{-1}\geq (q''-q')R^{2},\quad\rho^{2}\geq
\eta(q''-q')R^{2}.
\end{equation}
On the other hand, $Q\subset C_{R}(R^{2},0)$, $s>R^{2}$,
 and $\rho< R$.

Moreover, 
 by construction,  $|\Gamma\cap Q|\geq q_{0}|Q|$ and
by Lemma \ref{lemma 12.20.2} on the set where
$|z-y|\leq\rho/2$
$$
I(z):=E\Big(\int_{0}^{\tau}I_{\Gamma}(s+t,z+x_{s+t}-x_{s})
\,dt\mid\cF_{s}\Big)
\geq \mu_{0}\rho^{2}\geq \mu_{0}\eta(q''-q')R^{2},
$$
 where $\tau$ is the first exit time of $(s+t,z+x_{s+t}-x_{s})$
from $C_{ \rho}(s,y)$. Now by Theorem \ref{theorem 12.7.2}
for $x\in B_{\kappa R}$
$$
E \int_{0}^{\tau'_{R}(x)\wedge(2R^{2})}I_{\Gamma}(t,x+x_{t})\,dt
\geq E I_{\tau'_{R}(x)>s,|x+x_{s}-y|\leq\rho/2}
I(x+x_{s})
$$
$$
\geq \mu_{0}\eta(q''-q')R^{2}P_{x}\big
(\tau'_{R}(x)>s,|x+x_{s}-y|\leq\kappa\rho\big)\geq N^{-1}
(\rho/R)^{\nu}
\mu_{0}\eta(q''-q')R^{2}.
$$
This proves \eqref{7,27,1}.

By combining the two cases (i) and (ii) we conclude that
$$
G_{R}(\Gamma,x)\geq \min\big(\mu_{0}, \nu_{0}\mu_{R}(q''),
\mu_{0}\nu\eta^{n}(q''-q')^{n}\big)R^{2},
\quad |x|\leq \kappa R,
$$
and the arbitrariness of $\Gamma$ allows us to conclude that
\begin{equation}
                                               \label{7,27,5}
\mu_{R}(q')\geq \min\big(\mu_{0}, \nu_{0}\mu_{R}(q''),
\mu_{0}\nu\eta^{n}(q''-q')^{n}\big),
\end{equation}
whenever \eqref{1.5.1} holds. Observe that
\eqref{7,27,5} is identical to (9.3.10) of \cite{Kr_18}
and by literally repeating what is in \cite{Kr_18},
just replacing $\xi$ there with our $\theta$,
we come to \eqref{2.12.1} for any $R$.
By Lemma \ref{lemma 2.11.1} the right-hand side
in \eqref{2.12.1} can be taken the
same for $R\leq\bar R$. The theorem is proved.

The following four results are derived from Theorem
\ref{theorem 12.21.1} in the same way as similar results are derived
from Theorem 4.1 of \cite{Kr_19_1}.

\begin{corollary}
                               \label{corollary 10.11.1}
For any 
$\kappa\in(0,1)$ there exists  
$N$, depending only on $\kappa,d,\delta,p_{0},\uR,\bar R$,
and $\bar b_{\bar R}$,  such that, for any $R\in(0,\bar R]$,
$x\in   B_{\kappa R}$, and closed set
$\Gamma\subset  C_{R}(R^{2},0)$, the probability that the process 
$(t,x+x_{t})$  
reaches $\Gamma$ before exiting from $C_{2R^{2},R}$
is greater than or equal to $N^{-1} (|\Gamma|/|C_{R}|)^{\mu-1/
(d_{0}+1)}$:
\begin{equation}
                               \label{10.2.10}
P (\tau_{\Gamma}(x) <\tau_{2R^{2},R}(x) )
\geq N^{-1} (|\Gamma|/|C_{R}|)^{\mu-1/(d_{0}+1)},
\end{equation}
where $\tau_{\Gamma}(x) $ is the first time $(t,x+x_{t})$
hits $\Gamma$, $\tau_{2R^{2},R}(x)$ is the first exit time of
$(t,x+x_{t})$ from $C_{2R^{2},R}$,
 $\mu=1/\gamma$, and $\gamma$ is taken from Theorem \ref{theorem 12.21.1}.
\end{corollary}

\begin{corollary}
                      \label{corollary 10.1.1}    
For any $R\in(0,\bar R]$, Borel nonnegative $f$ 
vanishing outside $C_{R}(R^{2},0)$, and $x\in B_{\kappa R}$
$$
\int_{C_{R}(R^{2},0)}f^{1/(2\mu)}(t,y)\,dydt\leq NR^{d+2-1/\mu}
\Big(E \int_{0}^{\tau_{2R^{2},R}(x) }f(t,x+x_{t})\,dt\Big)^{1/(2\mu)},
$$
where $N$ depends only on $\kappa,d,\delta,p_{0},\uR,\bar R$,
and $\bar b_{\bar R}$,

\end{corollary}

\begin{corollary}
                                       \label{corollary 10.21.1}
 
Let $R\in(0,\bar R]$, $\gamma\in(0,1)$,
 and assume that a closed set $\Gamma\subset
B_{R}$ is such that, for any $r\in(0,R)$, $|B{r}\cap\Gamma|\geq
\gamma |B_{r}|$. Then there exist    constants $\alpha\in(0,1)$
and $N$, depending only on $\kappa,d,\delta,p_{0},\uR,\bar R$,
 $\bar b_{\bar R}$, and $\gamma$, such that, for any $ x\in B_{R/2}$,
\begin{equation}
                                     \label{10.21.1}
 P(\tau_{R}(x)<\tau_{\Gamma}(x))\leq N( |x| /R)^{\alpha},
\end{equation}
where $\tau_{\Gamma}(x)$ is the first time $ x+x_{t}$
hits $\Gamma$.
\end{corollary}
 
The fourth result has the same spirit
as Theorem 4.11 of \cite{Kr_19_1}
and can be used in investigating the boundary
behavior of solutions
of parabolic equations with drift in $L_{p,q}$.

We are going to use the following condition
\begin{equation}
                                                   \label{10.7.1}
p ,q \in [1,\infty],\quad
\nu:=1-\frac{d_{0}}{p }-\frac{1}{q }\geq 0,
\end{equation}
where $d_{0} \in(d/2,d)$,
depending only on
$\delta$, $d$, $\uR$,  $p_{0}$,  
  is taken from \cite{Kr_21_1}.   

 \begin{theorem}
                                        \label{theorem 10.20.3}
Let \eqref{10.7.1} be satisfied with $\nu=0$, $T\in(0,\infty)$,
 and
let $D$ be a bounded domain in $\bR^{d}$ with $0\in\partial D$. 
Assume that for some constants $\rho,\gamma>0$ and 
any $r\in (0,\rho)$ we have $|B_{r}\cap D^{c}|\geq \gamma
|B_{r}|$. Then there exist  $\beta >0$
and $N$, depending only on $ d,\delta,p_{0},\uR $,
  $\bar b_{\infty}$, $\gamma$ with 
 $N$ also depending on $\rho$ and the diameter of $Q:=(0,T)\times D$,
such that, for any nonnegative $f\in L_{p,q}(Q)$,
\begin{equation}
                                                 \label{10.20.4}
u(x):=E\int_{0}^{\tau(x)}f(t,x+x_{t})\,dt
\leq N|x|^{\beta}\|f\|_{L_{p,q}(Q)},
\end{equation}
where $\tau(x)$ is the first exit time of $(t,x+x_{t})$
from $Q$.
\end{theorem}
 
  In the next section we will need the following
fact of crucial importance, the origin
of which lies in \cite{KS_80} and \cite{Sa_80}.
A few other related results below also have their origin
 in \cite{KS_80} and
\cite{Sa_80} where the drift is bounded.

\begin{theorem}
                                      \label{thm:4.1.10} 
Let $\kappa,\eta,\zeta, q\in(0,1)$,  $R\in(0,\bar R]$,
$T\in[\eta R^{2},\eta^{-1}R^{2}]$,
and closed $\Gamma\subset C_{T,R}$ be such that
$|\Gamma\cap C_{\zeta T,R}((1-\zeta)T,0)|\geq q
|C_{\zeta T,R}|$. Then there exists
$\pi_{0}=\pi_{0}(\kappa,\eta,\zeta,q,d,\delta,p_{0},\uR,\bar R,
\bar b_{\bar R})>0$,
 such that,
for $(t ,x )\in C_{(1-\zeta)T, \kappa R}$,
\begin{equation}
                                      \label{12.27.3}
P (\tau_{\Gamma}(t ,x )<\tau_{T,R}(t ,x ))\geq \pi_{0},
\end{equation}
where $\tau_{\Gamma}(t ,x )$ is the
 first time $(t+s,x +x_{s})$
hits $\Gamma$ and $\tau_{T,R}(t ,x )$ is its
first exit time from $C_{T,R}$.
\end{theorem}

Proof.  
Observe that one can choose 
$\rho\in(0,1]$, depending only on $d,\eta,\zeta$, and $q$, and one can find
$(t^{0},x^{0})\in C_{T,R}$ with $t^{0}\ge \rho^{2}R^{2}
+(1-\zeta)T$
such that $C_{\rho R}(t^{0}+\rho^{2}R^{2},x^{0})\subset C_{T,R}$
and
 $|\Gamma\cap C_{\rho R}(t^{0}+\rho^{2}R^{2},x^{0})|
\geq \bar q|C_{\rho R} |$, where
$\bar q>0$ depends only on $d,\eta,\zeta$, and $q$. Then
by Corollary \ref{corollary 10.11.1}, for 
$x\in B_{  \kappa R}(x_{0})$
 the  probability that the process
$(t_{0} +s,x+x_{s})$ will hit $\Gamma$ before exiting
from $C_{2\rho^{2}R^{2},\rho R}(t^{0} ,x^{0})$ is estimated from 
below by a strictly positive constant depending only
on $\kappa,\bar q,d,\delta,p_{0},\uR,\bar R$, and $\bar b_{\bar R}$.
After that it only remains to invoke Theorem
\ref{theorem 12.7.2} recalling that $t^{0}\geq \rho^{2}R^{2}
+(1-\zeta)T$ and $t<(1-\zeta)T$.
The theorem is proved.

  \mysection{The case of diffusion processes}
                                   \label{section 10.24.1}

In this section, among other things, we generalize
some recent results in \cite{Ya_20} and extend them
to processes with singular drift.

Fix a constant $\delta\in(0,1)$ and recall that by
$\bS_{\delta}$ we denote the set of $d\times d$-symmetric
matrices whose eigenvalues are between
$\delta$ and $\delta^{-1}$. In this section
we impose the following.
\begin{assumption}
                                         \label{assumption 10.12.1}

(i) On $\bR^{d+1}$
 we are given Borel measurable $\sigma(t,x)$
and $b(t,x)$ with values  in $\bS_{\delta}$
and in $\bR^{d}$ respectively.

(ii) We are given $p_{0},q_{0} \in[1,\infty)$  
satisfying \eqref{5.10.1} and a function $h(t,x)$
satisfying \eqref{8.19.1} 
and such that $|b|\leq h$.

(iii) Assumption \ref{assumption 12.18.2}
is satisfied.
\end{assumption}

Let $\Omega$ be the set of $\bR^{d+1}$-valued
 continuous function $(t_{0}+t,x_{t})$, $t_{0}\in \bR$,
defined for $t\in[0,\infty)$.
For $\omega=\{(t_{0}+t,x_{t}),t\geq0 \}$, define
$\sft_{t}(\omega)=t_{0}+t$, $x_{t}(\omega)=x_{t}$,
and set $\cN_{t}=\sigma((\sft_{s},x_{s}),s\leq t)$,
$\cN=\cN_{\infty}$.  
In the following theorem which is Theorem
6.1 of \cite{Kr_20_2} we use the terminology from
\cite{Dy_63}.

\begin{theorem}
                                            \label{theorem 4.27.1}
 On $\bR^{d+1}$ there exists a strong Markov process
$$
X=\{(\sft_{t},x_{t}),\infty,\cN_{t}, P_{t,x})
$$
such that the process 
$$
X_{1}=\{(\sft_{t},x_{t}),\infty,
\cN_{t+}, P_{t,x})
$$
 is Markov and for any $(t,x)\in\bR^{d+1}$
there exists a $d$-dimensional Wiener process $w_{t}$, $t\geq0$,
which is a Wiener process relative to $\bar \cN_{t}$,
where $\bar \cN_{t}$ is the completion of $\cN_{t}$
with respect to $P_{t,x}$, and such that with 
$P_{t,x}$-probability one, for
all $s\geq 0$, $\sft_{s}=t+s$ and 
\begin{equation}
                                               \label{10.16.3}
x_{s}=x +\int_{0}^{s}\sigma(t  +r,x_{r})\,dw_{r}
+\int_{0}^{s}b(t  +r,x_{r})\,dr.
\end{equation}
 \end{theorem}

\begin{remark}
                                      \label{remark 10.29.1}

 To be completely rigorous, to refer to \cite{Kr_20_2}
we should have $b\in L_{p_{0},q_{0}}$ (globally),
and   \eqref{8.19.1} is not needed.
But with our $b$, owing to Corollary 2.8 and
Theorems \ref{theorem 9.7.1}, the arguments in \cite{Kr_20_2}
only simplify and do not require $b\in L_{p_{0},q_{0}}$.
Still it is worth saying that the author believes that
under only conditions in \cite{Kr_20_2} Harnack's inequality
is true. Regarding the H\"older continuity of caloric
functions in the same setting we have no guesses.
The H\"older continuity seems to require some sort
of self-similarity and the $ L_{p_{0},q_{0}}$-norm
is not preserved under such transformations
if $p_{0},q_{0}$ are subject to \eqref{5.10.1}.
\end{remark}

\begin{theorem}
                                 \label{theorem 9.24.1}
For any $\lambda\geq 1$, $ p,q $ satisfying 
\eqref{10.7.1}, and Borel nonnegative
$f(t,x)$  and for 
$$
 R_{\lambda}f(t,x):=E_{t,x}\int_{0}^{\infty}
e^{-\lambda s}f(t+s,x_{s})\,ds.
$$
we have
\begin{equation}
                                         \label{10.11.1}
\| R_{\lambda}f \|_{L_{p,q}(\bR^{d+1}_{+})}\leq N
\lambda^{-1} \| f \|_{L_{p,q}(\bR^{d+1}_{+})},
\end{equation}
where $N =N(\delta,d,\uR,p,q,p_{0},\bar b_{\infty} )$.

\end{theorem}

Proof.  
By Theorem \ref{theorem 8.30.1}
we have
\begin{equation}
                                         \label{10.23.1}
R_{\lambda}f (t,x)\leq 
N \lambda^{-\eta}
\|\Psi_{\lambda}^{1-\nu} f (t+\cdot,x+\cdot)\|_{L_{p ,q }
(\bR^{d+1}_{+})},
\end{equation}
where $\eta= \nu+( 2d_{0}-d)/(2p )$ and $\Psi_{\lambda}(t,x)
=\exp(- \sqrt{\lambda} 
(|x|+ \sqrt t)\bar\xi/16)$.
The right-hand side here coincides with the right-hand side
in (5.5) of \cite{Kr_21_1} and by borrowing the result
of  arguments from there
we get \eqref{10.11.1}. The theorem is proved.

 \begin{definition}
                                \label{definition 12.27.1}
If $Q$ is a set in $\bR^{d+1}$
and $u$ is a bounded Borel function on $Q$,
we call it caloric (relative to the process $X$) if
for any $(s,y)$ and $T,R\in(0,\infty)$
such that   $\bar C_{T,R}(s,y)\subset Q$
and any $(t ,x )\in C:=C_{T,R}(s,y)$ we have
$$
u(t ,x )=E_{t ,x } u(t +\tau_{C},x_{\tau_{C }}),
$$
where $\tau_{C}$ is the first exit time of $(t +
s,x_{s})$
from $C$.

\end{definition}

 First,
we deal with H\"older norm estimates for 
harmonic functions and potentials.
If $z_{1}=(t_{1},  x_{1})$ and $z_{2}=(t_{2}, x_{2})$,
we define    
\begin{equation}
                                                   \label{eq:4.2.5}
\rho(z_{1}, z_{2})=|x_{1}-x_{2}|+|t_{1}-t_{2}|^{1/2}
\end{equation}
 and call  $\rho(z_{1},z_{2})$   the parabolic
distance between $z_{1}$ and $z_{2}$.

\begin{lemma}
                                        \label{lem:4.2.2} 
Let $R\in(0, \bar R]$ and let $u$ be a caloric
function in $\bar C_{2R}$. Then there
exist constants $N$ 
and   
$$
\alpha_{0}\in(0,1),
$$ 
depending only on $\delta,d,p_{0},\uR,\bar R,
\bar b_{\bar R}$,
  such that, for any $\alpha\in(0,\alpha_{0}]$ and 
$z_{1},z_{2}\in C_{ R}$, we have   
\begin{equation}
                                       \label{eq:4.2.6}
\big|u(z_{1})-u(z_{2})\big|\le 
NR^{-\alpha}\rho^{\alpha}(z_{1},z_{2})
\sup\big(|u|,\bar C_{ 2R}\big).
\end{equation}

Furthermore, $\sup(|u|,\bar C_{ 2R})$ in 
\eqref{eq:4.2.6} can be replaced by $\osc(u,\bar 
C_{ 2R})$, where we use the 
notation
$$
 \osc(g,\Gamma) =\osc_{\Gamma}g=\sup_{\Gamma}g-
\inf_{\Gamma}g.
$$

\end{lemma}

Proof.  For $r$ such that $C_{r}\subset C_{2R}$, set
$$
w(r)=\osc(u,\bar C_{r}),\quad 
m(r)=\inf_{\bar C_{r}}u,\quad M(r)=\sup_{\bar C_{r}}u,
$$
$$
\mu(r)=(1/2)
\big(m(r)+M(r)\big).
$$

  Take  $r\le R/2$ and  suppose that
$$
\big|C_{ 2r}\cap\big\{ u\le
\mu(r)\big\} 
\big|\ge (1/2)|C_{ 2r}|.  
$$
Then there is a closed $\Gamma\subset
C_{ 2r}\cap\big\{ u\le
\mu(r)\big\}$ such that
\begin{equation}
                                   \label{eq:4.2.7}
\big|C_{3 r^{2},r}(r^{2},0)\cap\Gamma 
\big|\ge (1/4)|C_{3 r^{2},r}|
\end{equation}

By Theorem \ref{thm:4.1.10} for any $(t ,x )
\in \bar C_{r}$ we have
$$
P _{t,x}(\tau_{\Gamma} <\tau_{2r} )\geq \pi_{0},
$$
where $\pi_{0}>0$ depends only on $\delta,d,p_{0},\uR,\bar R,
\bar b_{\bar R}$,  
$\tau_{\Gamma} $ is the first time $(t  +s,x_{s})$
hits $\Gamma$, $\tau_{2r}$ is its first exit
time from $C_{2r}$.
Then by definition and the strong Markov property
for $\tau=\tau_{\Gamma}\wedge\tau_{2r}$ we have
$$
u(t ,x )= 
E_{t,x }  u(t +\tau_{2r} ,x_{\tau_{2r} })
$$
$$
=E_{t,x}  u(t+\tau _{2r},x_{\tau_{2r} })I_{\tau_{\Gamma}<\tau_{2r}}
+E_{t,x}  u(t+\tau _{2r},x_{\tau _{2r}})I_{\tau_{\Gamma}\geq\tau_{2r}}
$$
$$
=E_{t,x}  u(t+\tau _{\Gamma},x_{\tau_{\Gamma} })I_{\tau_{\Gamma}<\tau_{2r}}
+E_{t,x}  u(t+\tau_{2r} ,x_{\tau_{2r} })I_{\tau_{\Gamma}\geq\tau_{2r}}
$$
$$
\leq \mu(r)\pi_{0}+M(2r)(1-\pi_{0})
$$
(we used that $\mu(r)\leq M(2r)$).
It follows  that     
$$
M(r)\le \pi_{0}\frac{1}{2}\big(m(r)
+M(r)\big)+(1-\pi_{0})M(2r),
$$
$$
\big(1-\frac{\pi_{0}}{2})M(r)\leq \frac{\pi_{0}}{2} m(r)+(1-\pi_{0})M(2r).
$$

Adding  to this   the obvious inequality 
$$
\big(\frac{\pi_{0}}{2}-1)m(r)\leq -\frac{\pi_{0}}{2} m(r)
+(\pi_{0}-1)m(2r),
$$
we  get  
\begin{equation} 
                                             \label{eq:4.2.8}
\big(1-\frac{\pi_{0}}{2}\big)w(r)\le(1-\pi_{0})w(2r),\quad
w(r)\le\varepsilon w (2r),
\end{equation} 
where $\varepsilon<1$,  $\varepsilon=
\varepsilon(\pi_{0})$.  We may, certainly, assume that 
 $\varepsilon>1/2$.

 We have proved (\ref{eq:4.2.8}) assuming 
that  (\ref{eq:4.2.7}) is true. However if (\ref{eq:4.2.7}) 
is false, then  
$-u$ satisfies an inequality similar to \eqref{eq:4.2.7}
and  this leads to
 (\ref{eq:4.2.8}) again.

\label{iteration}
 Therefore,  $w(r)\le\varepsilon w(2r)$ for all $r\le R/2$.  
Iterations then yield    
$$
w(r)\le\varepsilon^{2}w(4r)\quad\text{for}\quad r\le R/4,..., 
w(r)\le\varepsilon^{n}w(2^{n}r)  \quad\text{for}\quad r\le2 ^{-n}R.
$$
 If $r\le R/2$ and  we take  $n:=\lfloor\log_{2}(R/r)\rfloor$,
 then  $r\le2^{-n}R$ and   
$$
w(r)\le\varepsilon^{n}w(2^{n}r)\le
\varepsilon^{-1}(r/R)^{\alpha}w(R)\le2
\varepsilon^{-1}(r/R)^{\alpha}
\sup\big(|u|,\bar C_{R }\big),
$$
where $\alpha=-\log_{2}\varepsilon\in(0,1)$.  
This provides an estimate of 
the oscillation of $u$ in any   $C_{r}$
with $r\le R/2$. The same  estimate
obviously holds for the oscillation of $u$ in  any   
  $ C_{ r}(t,x)\subset C_{ 2R}$ as long as  $r\le R/2$
and $(t,x)\in C_{R}$.

Now  take  $z_{1}=(t_{1},x_{1}),z_{2}=(t_{2},x_{2})\in C_{R}$
 such that
$r:=\rho(z_{1},z_{2})\le R/2$ and define 
$$
t=t_{1}\wedge t_{2},\quad x= (x_{1}+x_{2})/2.
  $$
 Then  we have $z_{i}\in\bar C_{ R}(t,x)$, $i=1,2$,
and     
\begin{align*}
\big|u(z_{1})-u(z_{2})\big|  \le &\, 2
\varepsilon^{-1}(r/R)^{\alpha}\sup\big(|u|,
 \bar C_{R }(t,x)\big)\\ 
   \le &\, 2\varepsilon^{-1}
\rho^{\alpha}(z_{1},z_{2})R^{-\alpha}\sup\big(|u|,\bar C_{ 2R}\big).
\end{align*}

 In  the case that $\rho(z_{1},z_{2})\geq R/2$  
 we have 
\begin{align*}
\big|u(z_{1})-u(z_{2})\big|   \le &\, 2
 \sup\big(|u|, \bar C_{2} \big)\\ 
  \le &\, 2^{1+\alpha}\rho^{\alpha}(z_{1},z_{2})R^{-\alpha}
\sup\big(|u|,\bar C_{ 2}\big).
\end{align*}

Thus,   $N=2^{1+\alpha}+2 
\varepsilon^{-1} $ in (\ref{eq:4.2.6})  is always
 a good choice
with   $\alpha $ found above. One can take any smaller
$\alpha$ as well since $\rho(z_{1},z_{2})\leq N(d)R$.
The lemma is proved.      

\begin{remark}
                                       \label{remark 2.12.1}
The constant $N$ in \eqref{eq:4.2.6}, generally,
depends on $\bar R$. However, if $\bar b_{\infty}
\leq \varepsilon$, where $\varepsilon>0$
depends only on $d$ and $\delta$, then this constant
is independent of $\bar R$. This is proved by using self-similar
transformations which change the process but  
allow us to take any $\uR$ we wish. In such situation
the Liouville theorem is valid: 
If $u$ is bounded and caloric in $\bR^{d+1}_{+} $,
then $u$ is constant (just send $R\to\infty$ in
\eqref{eq:4.2.6}).

\end{remark}

Here is the statement of the Harnack inequality.

\begin{theorem}
                                              \label{thm:4.2.1} 
Let   $R\in(0,\bar R]$, and let 
$u$ be a nonnegative 
caloric function in $\bar C_{2 R^{2},R}$.
 Then there exists a constant $N$, which depends 
only on $ \delta,d,\uR, \bar R, p_{0},\bar b_{\bar R}$, such that   
$$
u(R^{2},0)\le Nu(0,x) 
$$
whenever $|x|\le R/2$. 
\end{theorem}

Proof. We basically repeat the proof of  
Theorem 6.1 in  \cite{Kr_20} 
 and, to exclude a trivial situation, additionally assume that 
$$
u(R^{2},0)>0.
$$
 
  For  
$\kappa=1/2,\eta=1/2$, we  take  $N$ and $\nu$ from
Theorem \ref{theorem 12.7.2}, call this $N$ $N_{1}$,
 and,  having in mind
    Theorem \ref{theorem 11.8.1},
find $\gamma\in(0,1)$    close to 1 and $\varepsilon>0$
close to zero,   for which     
\begin{equation}
                                       \label{eq:4.2.3}
1-\varepsilon\geq q(\gamma)2^{-1}+\big[1-q(\gamma)  \big]2^{\nu}.
\end{equation}

  Next, for $r\in[0,R)$, introduce  
$$\mu(r)=u(R^{2},\,0)(1-r/R)^{-\nu},\quad 
n(r)=\sup\{ u,\bar C_{ r}(R^{2},0)\}
\,\, (n(0)=u(R^{2},0)),
$$
  and define
  $r_{0}$ as the greatest number in $r\in [0,R)$ satisfying 
$$
n(r)= \mu(r).
$$ 
 Such a number does exist because
  $n(0)=\mu(0)$, $\mu(r)\to\infty$
as $r\uparrow R$, and $n(r)$ is bounded, increasing,
and (H\"older) continuous.  
Choose $(t^{0},x^{ 0 })
\in \bar C_{ r_{0}}(R^{2},0)$
such that $n(r_{0})=
u(t^{ 0},x^{0})$ and consider  
the cylinder  
$$
Q:=\Big\{ (t,x)\,:\,0\leq t-t^{ 0}
<\frac{(R-r_{0})^{2}}{4},
\quad|x-x^{ 0}|<\frac{R-r_{0}}{2}\Big\} .
$$

  As is easy to see   $\bar Q\subset \bar C_{ r_{1}}(R^{2},0)$,
where $r_{1}=(R+r_{0})/2 $.  By the definition of $r_{0}$,
this implies   that    
$$
\sup_{\bar Q}u<\mu(r_{1})=u(R^{2},0)\Big(\frac{R-r_{0}}{2R}\Big)^{-\nu}
\leq 2^{\nu}n(r_{0}).
$$

We claim that  owing to this and 
 (\ref{eq:4.2.3}), 
\begin{equation}
                                                       \label{eq:4.2.4}
\big|Q\cap\big\{ u>n(r_{0})/2\big\}\big|\ge  (1- \gamma ) |Q|.
\end{equation}

 To argue by contradiction, assume \eqref{eq:4.2.4}
is false. Then
$$
\big|Q\cap\big\{ u\leq n(r_{0})/2\big\}\big|
>    \gamma   |Q|
$$
and there is a closed set $\Gamma\subset
Q\cap\big\{ u\leq n(r_{0})/2\big\}$ such that
$|\Gamma|> \gamma  |Q|$. 
Introduce $\tau_{\Gamma}$ as the first time the process
$(t^{0}+s,x_{t})$ hits $\Gamma$
and $\tau_{Q}$ as the first time it exits from $Q$.
It follows by definition, the strong Markov property
as in the proof of Lemma \ref{lem:4.2.2}, and from 
Theorem \ref{theorem 11.8.1}  that
(note that $n(r_{0})/2\leq\sup_{\bar Q}u$)
$$
u(t^{0},x^{0})=E_{t^{0},x^{0}}
I_{\tau_{\Gamma}<
\tau_{Q}}u(t^{0}+\tau_{\Gamma},
x_{\tau_{\Gamma}})
+E_{t^{0},x^{0}}I_{\tau_{\Gamma}
\geq\tau_{Q}}u(t^{0}+\tau_{Q},
x_{\tau_{Q}})
$$
$$
\leq P_{t^{0},x^{0}}
(\tau_{\Gamma}
<\tau_{Q})n(r_{0})/2+
(1-P_{t^{0},x^{0}}
(\tau_{\Gamma}
<\tau_{Q}))\sup_{\bar Q}u
$$
$$
\leq q(\gamma)n(r_{0})/2+
(1-q(\gamma))\sup_{\bar Q}u
$$
$$
\leq q(\gamma)n(r_{0})/2+
(1-q(\gamma))2^{\nu}n(r_{0}).
$$
Owing to \eqref{eq:4.2.3}
we  now have
$$
n(r_{0})\leq(1+\varepsilon)n(r_{0})
\big[q(\gamma)2^{-1}+(1-q(\gamma))2^{\nu}\big]
\leq (1-\varepsilon^{2})  n(r_{0}),
$$
which is impossible. This proves 
 (\ref{eq:4.2.4}).

  Next we  apply Theorem \ref{thm:4.1.10}
and get that 
$$
u(t^{0},x)\ge \pi_{0}n(r_{0})2^{-1}
$$
if $|x-x^{0}|\le (R-r_{0})4^{-1}$, where  
$\pi_{0}=\pi_{0}( d,\delta,p_{0},\uR,\bar R,\bar b_{\bar R},
 \gamma )>0$.   
 After that it only remains to apply  
Theorem \ref{theorem 12.7.2}
to conclude that for $|x|\le R$ we have    
$$
u(0,x)\ge\frac{1}{2}\pi_{0}n(r_{0})N_{1}^{-1}
\Big(\frac{R-r_{0}}{4}\Big)^{\nu}= 2^{-2\nu-1}\pi_{0}
N_{1}^{-1} u(4,0).   
$$
 The theorem is proved.

By using   Lemma \ref{lem:4.2.2} and Theorem \ref{theorem 9.5.1}
one derives in three lines the following analog
of Theorem 6.5 of \cite{Kr_20}.
 \begin{theorem}
                                     \label{theorem 10.8.1}
Assume that \eqref{10.7.1} holds with $\nu=0$. 
Let $R\in(0,\bar R/2]$ and
 let $g$ be a Borel bounded
function on $\bar C_{2R}$ and $f\in L_{p,q}(C_{2R})$.
For $(t ,x )\in C_{2R}$ define
\begin{equation}
                                                  \label{10.19.2}
u(t ,x )=E_{t ,x }
\int_{0}^{\tau_{2R}}f(t +s,x_{s})\,ds+
E_{t  ,x }g(t +\tau_{2R},x_{\tau_{2R}}),
\end{equation}
where $\tau_{2R} $ is the first exit time of 
$(t +s,x_{s})$ from $C_{2R}$.
  Then there exists a constant $N$, which depends
only on $\delta,d,\uR,\bar R, p, p_{0}$, and $\bar b_{\infty}$, such that   
\begin{equation}
                                                  \label{eq:4.2.9}
\big|u(z_{1})-u(z_{2})\big|\le N\big(R^{-\alpha}
\rho^{\alpha}(z_{1},z_{2})\sup_{\bar C_{ 2R}}|g|
+R^{(2d_{0}-d)/p}\|f\|_{L_{p,q}(C_{2R})}\big)
 \end{equation}
for $z_{1}$, $z_{2}\in  C_{ R}$,
$\alpha\in(0,\alpha_{0}]$, and $\alpha_{0}$
is taken from Lemma \ref{lem:4.2.2}.
 
\end{theorem}   

By playing with $R$ for fixed $z_{1}$, $z_{2}\in  C_{ R}$
as in the proof of Theorem 6.5 of \cite{Kr_20}
we get the following.

\begin{theorem}
                                    \label{theorem 12,14.2}
Under the conditions and notation from
 Theorem \ref{theorem 10.8.1}
there exists a constant $N$, which depends
only on $\delta,d,\uR,\bar R, p, p_{0}$, and $\bar b_{\bar R}$,
 such that    
\begin{equation}
                                           \label{12,14.6}
\big|u(z_{1})-u(z_{2})\big|\le NR^{-\beta}
\rho^{\beta}(z_{1},z_{2})
\big(\sup_{\bar C_{ 2R}}|u|
+R^{(2d_{0}-d)/p}\|f\|_{L_{p,q}(C_{2R})}\big)
\end{equation}
for $z_{1}$, $z_{2}\in  C_{R}$, where
$$
\beta=\frac{\alpha_{0}(2d_{0}-d)}{\alpha_{0} p+2d_{0}-d}.
$$

\end{theorem}

As a standard consequence of just continuity of $u$
we have the following.
\begin{theorem}
                                    \label{theorem 10.28.1}
The process
$$
X_{1}=\{(\sft_{t},x_{t}),\infty,
\cN_{t+}, P_{t,x})
$$
is strong Markov.
\end{theorem}

\mysection{Applications}
                                          \label{section 2.17.1}

Here we suppose that Assumption \ref{assumption 10.12.1}
is satisfied and set $a=\sigma^{2}$,
$$
Lu(t,x)=(1/2)a^{ij}(t,x)D_{ij}u(t,x)+b^{i}(t,x)D_{i}u(t,x).
$$
\begin{theorem}
                                        \label{theorem 10.19.1}
Let   $R\in(0,\bar R/2]$ and
assume that \eqref{10.7.1} holds with $\nu=0$, $p<\infty$,
$q<\infty$  and that
we are given a function 
$u\in W^{1,2}_{p,q,\loc}(C_{2R})\cap C(\bar C_{2R})$.
  Then for $-f=\partial_{t}u+Lu$ we have  
\begin{equation}
                                                  \label{10.19.1}
\big|u(z_{1})-u(z_{2})\big|\le NR^{-\beta}
\rho^{\beta}(z_{1},z_{2})
\big(\sup_{\bar C_{ 2R}}|u|
+R^{(2d_{0}-d)/p}\|f\|_{L_{p,q}(C_{2R})}\big)
 \end{equation}
for $z_{1}$, $z_{2}\in  C_{ R}$, where $N$ and $\beta$ 
are taken from
Theorem \ref{theorem 12,14.2}.
\end{theorem}

Proof. Approximating $C_{2R}$ by $C_{2R-\varepsilon}$
we see that we may assume that 
$u\in W^{1,2}_{p,q }(C_{R})\cap C(\bar C_{R})$.
  This gives us the opportunity
to replace $L$ in  the definition of $f$
with $L_{n}:=I_{|b|\geq n}\Delta+I_{|b|< n}L$
and then pass to the limit
by the dominated convergence and monotone convergence theorems.
Hence, we may assume that  $b$ is bounded.
After that it only remains to use 
It\^o's formula (Theorem \ref{theorem 10.15.1})
 for the Markov process from
Section \ref{section 10.24.1} (cf.~\eqref{10.16.3})
to see that $u$ has form \eqref{10.19.2}
for which \eqref{12,14.6} is valid. The theorem is proved.

\begin{remark}
                                        \label{remark 10.19.1}
A consequence of Theorem \ref{theorem 10.19.1}
is a rather weak statement that any
$u\in W^{1,2}_{p,q,\loc}(C_{R})$ admits a modification
which is in $C^{\beta}_{\loc}(C_{R})$.

Indeed, if  $u\in W^{1,2}_{p,q,\loc }(C_{ R})$
then its mollifiers will belong to $W^{1,2}_{p,q  }
(C_{ R-\varepsilon})$ and by \eqref{12,14.6} 
{\em with\/} $L=\Delta$ will be  in $C^{\beta} (C_{
(R-\varepsilon)/2})$. Passing from the mollifiers
to the function itself we find the modification
in question in $C_{R/2}$. After that scaling and
shifting the origin takes care of the rest of $C_{R}$.

\end{remark}

\begin{remark}
                                        \label{remark 10.24.1}
If $u$ is bounded in $C_{2R}$, belongs to
$ W^{1,2}_{p,q,\loc}(C_{2R})$, and is caloric
($\partial_{t}u+Lu=0$) in $C_{2R}$, then Theorem
\ref{theorem 10.19.1} implies that it is
H\"older continuous in $C_{R}$ with the exponent
and constant independent of any regularity
of $a$ and $b$. This fact along with Harnack's inequality
 was first proved
in \cite{KS_80} for bounded $b$ and $p=q=d+1$
in the parabolic case and $p=d$ in the elliptic case. They
were 
generalized by Cabr\'e \cite{Ca_95},
 Escauriaza \cite{Es_93}, and Fok \cite{Fo_95} in the elliptic
case when $p<d$ (close to $d$) again
  when $b$ is   bounded. In \cite{CKS_00}
parabolic version of these results, extending some earlier
results by Wang, are given
for $L_{p}$-viscosity solutions with
$p<d+1$ (close to $d+1$) when $b$ is   bounded.
 In our situation
we have some freedom in choosing $p,q$ and $b\in L_{p_{0},q_{0}}$,
but we only treat true solutions.

\end{remark}
\begin{theorem}
                                        \label{theorem 10.19.3}
Let   $R\in(0,\bar R]$ and
assume that \eqref{10.7.1} holds with $\nu=0$, $p<\infty$,
$q<\infty$.
Let   
$u\in W^{1,2}_{p,q,\loc}(C_{2 R^{2},R}) \cap 
C(\bar C_{2 R^{2},R})$
be such that $u>0$ on $\partial'C_{2 R^{2},R}$.
 Then there exists a constant $N$, which depends 
only on $ \delta,d,\uR,\bar R, p, p_{0}$,
and $\bar b_{\bar R}$, such that   
$$
u(R^{2},0)\le Nu(0,x)+N
R^{(2d_{0}-d)/p}\|f\|_{L_{p,q}(C_{2 R^{2},R})}
$$
whenever $|x|\le R/2$, where $-f=\partial_{t}u+Lu$.  
In particular, if  $\partial_{t}u+Lu=0$ in
$C_{  2R^{2},R}$ (a.e.), then (Harnack's inequality) 
$$
u(R^{2},0)\le Nu(0,x).
$$ 
\end{theorem}

Proof. As in the proof of Theorem \ref{theorem 10.19.1},
the general case is reduced to the one in which 
$b$ is bounded and 
$u\in W^{1,2}_{p,q }(C_{2 R^{2},R})
\cap C(\bar C_{2 R^{2},R} )$. In that case,
as in the proof of Theorem \ref{theorem 10.19.1},
by It\^o's formula for the Markov process from
Section \ref{section 10.24.1} 
\begin{equation}
                                              \label{11.3.2}
u(t,x)=h(t,x)+F(t,x),
\end{equation}
where 
$$
h(t,x)=E_{t,x}u(t+\tau,x_{\tau})\geq0,\quad
F(t,x)=E_{t,x}\int_{0}^{\tau}f(t+s,x_{s})\,ds,
$$
and $\tau  $ is the first exit time of 
$(t +s,x_{s})$ from $C_{2 R^{2},R}$.
By Theorem \ref{thm:4.2.1}, $h(R^{2},0)\le Nh(0,x)$
and it only remains to use Theorem \ref{theorem 9.5.1}
to estimate $F$. The theorem is proved.

Here is a generalization
of the Fanghua Lin estimate for operators
with summable drift which is one of the main tools in
the Sobolev space theory of fully nonlinear parabolic
equations (see, for instance, \cite{Kr_18}). 

\begin{theorem}
                       \label{theorem 10.6.10}

Let  $R\in(0,\bar R]$,  $p,q$
satisfy \eqref{10.7.1} with $\nu=0$, $p<\infty$,
$q<\infty$. Let
$u\in W^{1,2}_{p,q,\loc}(C_{R})\cap C(\bar C_{R})$, and
$c\in L_{p,q}(C_{R})$,
$c\geq0$. Then  
$$
\Big(\dashint_{C_{R}}\big(|D^{2}u|+(|b|+R^{-1}) |Du|\big)^{1/(2\mu)} 
\,dxdt\Big)^{2\mu}
$$ 
\begin{equation}
                              \label{10.6.1}
 \leq NR^{-d/p-2/q} \|f\|_{L_{p,q}(C_{R})}
 +NR^{ -2}\sup_{\partial' C_{R}}|u|,
\end{equation}
where $-f=\partial_{t}u+Lu-cu$,
$\mu$ is taken from Corollary \ref{corollary 10.1.1} 
with $\kappa=1/2$
and $N$ depends only on
$ d,\delta,\uR, \bar R, p, p_{0} $, 
  $R^{2-d/p-2/q}\|c\|_{L_{p,q}(C_{R} )}$,  $\bar b_{\uR}$,
$\bar b_{\bar R}$,  and
the function  $\bar N(d,p_{0},\cdot)$
(see \eqref{12.18.3}).
\end{theorem}

Proof.  
On the account of moving $R$,
we may assume that 
$u\in W^{1,2}_{p,q}(C_{R})$.
After that we observe that in light of 
Theorem \ref{theorem 10.14.1}
$$
\|\partial_{t}u+Lu\|_{L_{p,q}(C_{R})}\leq\|
\partial_{t}u+Lu-cu\|_{L_{p,q}(C_{R})}
+\| c\|_{L_{p,q}(C_{R})}\sup_{C_{R}}|u|
$$
$$
\leq \big(1+NR^{2-d/p-2/q}\| c\|_{L_{p,q}(C_{R})}\big)
\|\partial_{t}u+Lu-cu\|_{L_{p,q}(C_{R})}
+\| c\|_{L_{p,q}(C_{R})}\sup_{\partial' C_{R}}|u|
$$
and reduce the case of general $c$ to the one with $c\equiv0$.
As a few times before we may assume that $b$ is bounded and then
using approximations we see that assuming that
$u\in C^{1,2}(\bar C_{R})$
and that the coefficients   $a^{ij}$ are infinitely differentiable
do not restrict generality. In that case introduce 
$L'(t,x)=L(t,x)$ for $t\geq 0$ and $L'(t,x)=L(-t,x)$
for $t<0$ and introduce
$v$
as a unique $W^{1,2}_{d+1}(C_{2R^{2},R}(-R^{2},0))$-solution
of the equation 
$$
\partial_{t}v+L'v=-fI_{C_{R}}
$$
with boundary condition $v=u$ on $\{t\geq 0\}
\cap\partial'C_{2R^{2},R}(-R^{2},0)$ and 
$v(t,x)=u(-t,x)$ on $\{t\leq 0\}
\cap\partial'C_{2R^{2},R}(-R^{2},0)$. Owing to uniqueness
$v=u$ in $C_{R}$ and by Theorem \ref{theorem 10.14.1}
in $C_{2R^{2},R}(-R^{2},0)$ we have
\begin{equation}
                                               \label{10.25.1}
|v|\leq NR^{2-d/p-2/q}
\|f\|_{L_{p,q}(C_{R})}+\sup_{\partial'C_{R}}|u|.
\end{equation}

Next, it is easy to see that for sufficiently
small $\varepsilon>0$, depending only on $\delta$
and the function $\bar N(d,p_{0},\cdot)$, we have
   that
$\hat a:= a-\varepsilon I_{C_{R}}(D_{ij}u)/|D^{2}u|\in 
\bS_{(\delta^{2}-\varepsilon)^{1/2}}$ and
\begin{equation}
                                                     \label{2.15.1}
\bar N(d,p_{0},(\delta^{2}-\varepsilon)^{1/2})\bar b_{\uR}
<1
\end{equation}
(see \eqref{12.18.3} and Assumption \ref{assumption 12.18.2}).
Furthermore, for $\hat b=b-\varepsilon I_{C_{\bar R}}(|b|+1)Du/|Du|$
any $\rho>0$ and $(t,x)\in\bR^{d+1}$ we have
$$
\|\hat b\|_{L_{p_{0},q_{0}}(C_{\rho}(t,x))}
\leq(1+\varepsilon)\| b\|_{L_{p_{0},q_{0}}(C_{\rho}(t,x))}
+\varepsilon N(d)(\rho\wedge\bar R).
$$
It follows that
$$
\|\hat b\|^{q_{0}}_{L_{p_{0},q_{0}}(C_{\rho}(t,x))}
\leq \Big((1+\varepsilon)\bar b_{\rho}^{1/q_{0}}+
N\varepsilon\Big)^{q_{0}}\rho,
$$
where $N$ depends only on $d,p_{0}$, and $\bar R$.
It is seen that for a $\varepsilon>0$,
depending only on $d,\delta,p_{0}$, $\bar R$,
and the function $\bar N(d,p_{0},\cdot)$, not only 
\eqref{2.15.1} is satisfied but also
$$
\bar N(d,p_{0},(\delta^{2}-\varepsilon)^{1/2})
\Big((1+\varepsilon)\bar b_{\uR}^{1/q_{0}}+
N\varepsilon\Big)^{q_{0}}<1.
$$

Therefore the above theory is applicable
to the operator
$$
\hat L =(1/2)\hat a^{ij}D_{ij}+\hat b^{i}D_{i}.
$$
 Then set
 $\hat \sigma
=\hat a^{1/2}$ and 
  consider the diffusion process $(-R^{2}+t,x_{t})$, $t\geq0$,
starting from $(-R^{2},0)$ with  diffusion matrix $\hat \sigma$
and drift  $\hat b$.
By It\^o's formula
$$
v(-R^{2},0)=-E\int_{0}^{\tau}(\partial_{t}v+\hat Lv)
(-R^{2}+t,x_{t})\,dt
+Ev(-R^{2}+\tau,x_{\tau}),
$$
where $\tau$ is the first exit time of
$(-R^{2}+t,x_{t})$ from $C_{2R^{2},R}(-R^{2},0)$. Obviously,
$$
|Ev(-R^{2}+\tau,x_{\tau})|\leq\sup_{\partial' C_{R }}|u|.
$$
Furthermore, on $C_{R }$ we have 
$$
\partial_{t}v+\hat Lv
=\partial_{t}u+\hat Lu=-f-\varepsilon|D^{2}u|-\varepsilon
(|b|+1)|Du|.
$$
In the remaining part of $C_{2R^{2},R}(-R^{2},0)$ we have
$\partial_{t}v+\bar Lv=0$. It follows that
$$
 E\int_{0}^{\tau}I_{C_{1}}\varepsilon\big(|D^{2}u|
+(|b|+1)|Du|\big)(-R^{2}+t,x_{t})\,dt
$$
$$
\leq v(-R^{2},0)-E\int_{0}^{\tau}I_{C_{R}}f(-R^{2}+t,x_{t})\,dt
+\sup_{\partial' C_{R}}|u|.
$$
After that it only remains to recall \eqref{10.25.1}
and apply Theorem   \ref{theorem 10.14.1} and
Corollary \ref{corollary 10.1.1}.
The theorem is proved.

\mysection{Appendix}
                                     \label{section 2.14.1}
 
Here we present without proofs some results
from \cite{Kr_21_1} frequently used in the main text.

Set 
\begin{equation}
                                                  \label{1.6.1}
\tau'_{R}(x)=\inf\{t\geq0:x+x_{t}\not\in B_{R}\},\quad
\gamma_{R}(x)=\inf\{t\geq0:x+x_{t} \in \bar B_{R}\}.
\end{equation}

\begin{theorem}[Theorem 2.3]
                                      \label{theorem 8.2.1}
There are
   constants $\bar \xi=\bar \xi(d,\delta)\in (0,1) $ and
 $\bar N=\bar N(d,p_{0},  \delta)$ {\em
continuously\/} depending on $\delta$
such that if, for an $  R\in(0,\infty)$, we have
\begin{equation}
                                     \label{12.18.3}
\bar N \bar b_{  R}\leq 1,
\end{equation}
then  for $|x|\leq R$
\begin{equation}
                                          \label{8.2.2} 
  P( \tau_{R}(x)  =   R^{2} )\leq 1-\bar\xi,
\quad P( \tau_{R}  =   R^{2} )\geq 
 \bar\xi .   
\end{equation}
Moreover for $n=1,2,...$ and $|x|\leq R$
\begin{equation}
                                          \label{1.3.1} 
P(\tau'_{R}(x)\geq nR^{2})
=  P( \tau_{nR^{2},R}(x)  =  n R^{2} )\leq (1-\bar\xi)^{n},   
\end{equation}
so that $E\tau'_{R}(x)\leq N(d,\delta)R^{2}$.

Furthermore,  for any $x\in \bar B_{9R/16}$
\begin{equation}
                                          \label{1.2.1} 
P(\tau'_{R}(x)>\gamma_{R/16}(x))\geq\bar\xi.   
\end{equation}

\end{theorem}

\begin{theorem}[Theorem 2.6]
                                     \label{theorem 8.20.1}
For any $\lambda,R>0$ we have
\begin{equation}
                                          \label{8.20.1}
Ee^{-\lambda  \tau_{R} }\leq
e^{\bar\xi/2}e^{- \sqrt{\nelambda}  R
\bar\xi/2}=
\begin{cases}
 e^{\bar\xi/2}e^{-\sqrt\lambda R\bar\xi/2}
\quad\text{if}\quad \lambda
 \geq \ulambda\\
e^{\bar\xi/2}e^{-\lambda R\uR \bar\xi/2}\quad\text{if}
\quad \lambda
 \leq \ulambda,
\end{cases}
\end{equation}
where 
$$
\nelambda= \lambda\min(1, \lambda/\ulambda ),\quad \ulambda=\uR^{-2}.
$$ 
In particular, for  
  any   $R>0$ and $t\leq R\uR
\bar \xi/4 $ we have
\begin{equation}
                                             \label{10.2.2}
P(  \tau_{R}  \leq  t )\leq 
 e^{\bar\xi/2}\exp\Big(-\frac{{\bar \xi}^{2}R^{2}}{16 t}\Big).
\end{equation}

\end{theorem}

\begin{theorem}[Theorem 2.9]
                                        \label{theorem 1.24.1}
Let $R\in(0,\uR]$, $x,y\in\bR^{d}$ and $16|x-y|\geq 3R$.  
For $r>0$ denote by $S_{r}(x,y)$   the open convex hull
of $B_{r}(x)\cup B_{r}(y)$. Then there exist
$T_{0},T_{1}$, depending only on $\bar\xi$,
such that $0<T_{0}<T_{1}<\infty$ and the probability $\pi$
that $x+x_{t}$ will reach $\bar B_{R/16}(y)$ before exiting
from $S_{R}(x,y)$ and this will happen
on the time interval $[nT_{0}R^{2},nT_{1}R^{2}]$
is greater than $\pi_{0}^{n}$, where
$$
n= \Big\lfloor \frac{16|x-y|+R}{4R}\Big\rfloor 
$$  
and $\pi_{0}=\bar\xi/3$.
\end{theorem}

\begin{theorem}[Theorem 4.3]
                                           \label{theorem 2.3.1}
There exists $d_{0}\in(1,d)$, depending only on
$\delta$, $d$, $\uR$,  $p_{0}$,  
 such that for any $p\geq d_{0}+1$ and
$\lambda>0$
$$
\int_{0}^{\infty}\int_{\bR^{d}}G_{\lambda}^{p/(p-1)}(t,x)
\,dxdt\leq N(\delta,d,\uR,p_{0},\lambda,p).
$$
Furthermore, the above constant
$N(\delta,d,\uR,p_{0},\lambda,p)$
can be taken in the form
$$
N(\delta,d,\uR,p_{0},p)\melambda_{p}^{(d+2)/(2p)-1},
$$
where
$$
\melambda_{p}=\lambda(1\wedge\lambda)^{d/(2p-d-2)}.
$$
\end{theorem}

\begin{theorem}[Theorem 4.8]
                                  \label{theorem 8.30.1}
Suppose   
\begin{equation}
                                                   \label{10.7.13}
p ,q \in [1,\infty],\quad
\nu:=1-\frac{d_{0}}{p }-\frac{1}{q }\geq 0.
\end{equation}
Then there is  $N=N(\delta,d,\uR,p,q,p_{0},\bar b_{\infty} )$ such that
for any $\lambda>0$ and   Borel nonnegative $f$  
we have 
\begin{equation}
                                                   \label{8.22.40}
 E\int_{0}^{\infty}e^{- \lambda t}
f(t,x_{t}) \,dt\leq 
N\melambda_{d_{0}+1}^{-\nu+(d-2d_{0})/(2p )}\|\Psi_{\lambda}^{1-\nu} f\|_
{L_{p ,q }(\bR^{d+1}_{+})},   
\end{equation}
where $\Psi _{\lambda}(t,x)=\exp(- 
\sqrt{\nelambda} (|x|+ \sqrt t)\bar\xi/16)$.
In particular, if $f$ is independent of $t$, $p\geq d_{0}$,
 and $q=\infty$
$$
E\int_{0}^{\infty}e^{- \lambda t}
f( x_{t}) \,dt\leq 
N\melambda_{d_{0}+1}^{-1+d/(2p)}\|\bar \Psi_{\lambda}^{d_{0}/p} f\|_
{L_{p   }(\bR^{d } )},
$$
where $\bar \Psi _{\lambda}( x)=\exp(- 
\sqrt{\nelambda}  |x| \bar\xi/16)$.  
\end{theorem}

\begin{theorem}[Theorem 4.9]
                                          \label{theorem 9.7.1}
Assume that \eqref{10.7.13} holds.
Then

(ii)
  for any
$n=1,2,...$, nonnegative Borel $f$ on $\bR^{d+1}_{+}$, and
 $T\leq 1$  we have
\begin{equation}
                                          \label{9.7.1}
E\Big[\int_{0}^{T}  
f(t,x_{t})\,dt\Big]^{n}\leq n!N^{n} 
T^{n\chi }\| \Psi^{(1-\nu)/n} _{1/T}
f\|^{n}_{L_{p,q}(\bR^{d+1}_{+}) },
\end{equation}
where $N=N(\delta,d,\uR,p,q,p_{0},\bar b_{\infty} )$ and 
$\chi=\nu+(2d_{0}-d)/(2p)$; 

(ii)   for any
  nonnegative Borel $f$ on $\bR^{d+1}_{+}$, and
 $T\geq 1$  we have
\begin{equation}
                                          \label{9.7.10}
I:=E \int_{0}^{T}  
f(t,x_{t})\,dt \leq N  T^{1-1/q}
 \| \Psi^{  1-\nu   } _{1}
f\| _{L_{p,q}(\bR^{d+1}_{+}) },
\end{equation}
where $N=N(\delta,d,\uR,p,q,p_{0},\bar b_{\infty} )$.
 
\end{theorem}

\begin{theorem}[Theorem 4.10]
                                       \label{theorem 9.5.1}
Assume that \eqref{10.7.13} holds with $\nu=0$.
 Then 
for any $R\in(0,\bar R] $, $x$,  and Borel nonnegative $f$ 
given on $C_{R}$,
we have
\begin{equation}
                                                     \label{9.5.4}
E\int_{0}^{\tau_{R}(x)}f(t,x+x_{t})\,dt\leq
NR^{(2d_{0}-d)/p}\|f\|_{L_{p,q}(C_{R})},
\end{equation}
where $N=N(\delta,d,\uR,p, p_{0},\bar b_{\bar R},\bar R )$.

\end{theorem}

\begin{theorem}[Theorem 4.11]
                                 \label{theorem 10.15.1}
Assume that \eqref{10.7.13} holds with $\nu=0$
and  $p<\infty$,
$q<\infty$. Let
$Q$ be a bounded domain in $\bR^{d+1}$, $0\in Q$,
$b$ be {\em bounded\/}, and  $u\in W^{1,2}_{p,q}(Q)\cap C(\bar Q)$. Then,
for $\tau$ defined as the first exit time of $(t,x_{t})$
from $Q$ and for all $t\geq0$,
$$
u(t\wedge\tau,x_{t\wedge\tau})
=u(0,0)+\int_{0}^{t\wedge\tau}D_{i}u(s,x_{s})\,dm^{i}_{s}
$$
\begin{equation}
                                      \label{10.15.1}
+\int_{0}^{t\wedge\tau}[
\partial_{t}u(s,x_{s})+ a^{ij}_{s}D_{ij}u(s,x_{s})
+b^{i}_{s}D_{i}u(s,x_{s})]\,ds
\end{equation}
and the stochastic integral above is a square-integrable
martingale.
\end{theorem}

\begin{theorem}[Theorem 5.1]
                                        \label{theorem 10.14.1}
Let $0<R\leq\bar R$, domain $Q\subset C_{R}$,  and
assume that \eqref{10.7.1} holds with $\nu=0$, $p<\infty$,
$q<\infty$, and that
we are given a function 
$u\in W^{1,2}_{p,q,\loc}(Q)\cap C(\bar Q)$.
  Take a function
$c\geq 0$ on $Q$. Then on $ Q$
\begin{equation}
                                               \label{10.14.10}
u \leq NR^{(2d_{0}-d)/p}
\|I_{Q,u>0}(\partial_{t}u+Lu-cu)_{-}\|_{L_{p,q} }
+\sup_{\partial'Q}u_{+},
\end{equation}
where $N=N(\delta,d,\uR,\bar R,p, p_{0},\bar b_{\bar R})$
and $\partial'Q$ is the parabolic boundary of $Q$.
In particular (the maximum principle),
if $\partial_{t}u+Lu-cu \geq0$ in $Q$ and $u\leq0$
on $\partial'Q$, then $u\leq 0$ in $Q$.
\end{theorem}

\begin{lemma}[Lemma 2.2]
                                      \label{lemma 8.16.1}
 
We have
\begin{equation}
                                          \label{8.16.1}
A:=
E \tau_{R} (x)\leq  R^{2},
\end{equation}
 and,  assuming that \eqref{10.7.13} holds with $\nu=0$,
 for any Borel nonnegative $f$ we have
\begin{equation}
                                          \label{9.29.2}
 E\int_{0}^{\tau_{R}(x)}  
f( t,x_{t})\,dt\leq 
 N(d,p_{0},  \delta)(1+\bar b_{R}) ^{d/  p  }
 R 
 ^{d/  p  }
\|f\|_{L_{p ,q }}. 
\end{equation}
\end{lemma}

\begin{corollary}[Corollary 2.10]
                                     \label{corollary 1.25.1}
Let $R\leq\uR$, $\kappa\in[0,1)$, and $|x|\leq\kappa R$.
Then for any $T>0$
\begin{equation}
                                        \label{1.25.2}
NP(\tau'_{R}(x)> T)\geq e^{-\nu T/[(1-\kappa)R]^{2}},
\end{equation}
where $N$ and $\nu>0$ depend only on $\bar\xi$.
\end{corollary}

\begin{corollary}[Corollary 2.7]
                                       \label{corollary 7.29.1}
Let   $  \Lambda\in(0,\infty)$.
Then there is a  constant   
  $N =N (\uR,\bar R,\Lambda,\bar\xi)$
such that for any $R\in(0,\bar R]$, $\lambda\in[0,\Lambda]$ 
\begin{equation}
                                                   \label{8.21.1}
N E \tau_{R} \geq R^{2},\quad N E\int_{0}^{\tau_{R}}
e^{-\lambda t} \,dt  \geq R^{2} .
\end{equation}
\end{corollary}

\end{document}